\documentclass[10pt]{article}

\usepackage{amsfonts,latexsym}

\newcommand{\R}{\mathbb R}

\newcommand{\del}{\partial}
\newcommand{\e}{\varepsilon}

\newtheorem{theorem}{Theorem}[section]
\newtheorem{lemma}{Lemma}[section]
\newtheorem{proposition}{Proposition}[section]

\newtheorem{definition}{Definition}[section]

\begin{document}

\title{Singular limits for $4$-dimensional semilinear elliptic problems with
exponential nonlinearity}

\date{}

\author{Sami Baraket, Makkia Dammak, Taieb Ouni and  Frank Pacard}

\maketitle

\begin{abstract}
Using some nonlinear domain decomposition method, we prove the
existence of singular limits for solution of semilinear elliptic
problems with exponential nonlinearity.
\end{abstract}

\section{Introduction and statement of the results }

In the last decade important work has been devoted to the
understanding of singularly perturbed problems, mostly in a
variational framework. In general, a Liapunov-Schmidt type reduction
argument is used to reduce the search of solutions of singularly
perturbed partial differential equations to the search of critical
points of some function which is defined over some finite
dimensional domain.

\medskip

One of the purpose of the present paper is to present a rather
efficient method to solve such singularly perturbed problems. This
method has already been used successfully in geometric
context(constant mean curvature surfaces, constant scalar curvature
metrics, extremal K\"ahler metrics, manifolds with special holonomy,
\ldots ) but has never appeared in the context of partial
differential equations. We felt that, given the interest in singular
perturbation problems, it was worth illustrating this on the
following problem~:

\medskip

Let $\Omega \subset {\mathbb{R}}^{4}$ be a regular bounded open
domain in ${\mathbb{R}}^{4}$. We are interested in positive
solutions of
\begin{equation} \left\{\begin{array}{rclll}\Delta^{2} u &=&
\rho^{4} e^{u} &\mbox{ in}&~ \Omega\\\\u& =&\Delta u =0 &\mbox{on}&~
\partial\Omega,
\end{array}
\right. \label{eq:1.1}
\end{equation}
when the parameter $\rho$ tends to $0$. Obviously, the application
of the implicit function theorem yields the existence of a smooth
one parameter family of solutions $(u_\rho)_\rho$ which converges
uniformly to $0$ as $\rho$ tends to $0$. This branch of solutions is
usually referred to as the branch of {\em minimal solutions} and
there is by now quite an important literature which is concerned
with the understanding this particular branch of solutions
\cite{Mignot-Puel}.

\medskip

The question we would like to study is concerned with the existence
of other branches of solutions as $\rho$ tends to $0$. To describe
our result, let us denote by  $G(x, \cdot)$ the solution of
\begin{equation}
\left\{\begin{array}{rllllllll} \Delta^{2} G(x, \cdot ) & = & 64
\, \pi^{2} \, \delta_{x}  &  & & \mbox{ in} \quad  \Omega \\[3mm]
G(x, \cdot )& = & \Delta \, G(x, \cdot) & = & 0 & \mbox{on} \quad
\partial\Omega.
\end{array}\right.\label{eq:1.2}
\end{equation}
It is easy to check that the function
\begin{equation}
R(x, y) := G(x, y) + 8 \, \log |x-y|
\label{eq:1.3}
\end{equation}
is a smooth function.

\medskip

We define
\begin{equation}
W (x^{1},\ldots ,x^m)  : = \sum_{j=1}^{m} R(x^{j},x^{j}) + \sum_{j
\neq \ell} G(x^{j},x^{\ell}). \label{eq:1.4}
\end{equation}

Our main result reads~:
\begin{theorem}
Assume that $(x^1,\ldots ,x^m)$ is a {\em nondegenerate} critical
point of $W$, then there exist $\rho_{0} > 0$ and $(u_{\rho})_{\rho
 \in (0 , \rho_{0}) }$ a one parameter family of solutions of (\ref{eq:1.1}),
such that
\[
\lim_{\rho\rightarrow 0} u_{\rho} = \sum_{j=1}^m G(x^j , \cdot)
\]
in ${\mathcal{C}}_{loc}^{4,\alpha}(\Omega - \{x^1 , \ldots ,x^m \}
\, )$. \label{th:1.1}
\end{theorem}
This result is in agreement with the result of Lin and Wei
\cite{Lin-Wei} (see also \cite{wei}) where sequence of solutions
of (\ref{eq:1.1}) which blow up are studied. Indeed, in this
paper, the authors show that blow up points can only occur at
critical points of $W$.

\medskip

Our result reduces the study of nontrivial branches of solutions of
(\ref{eq:1.1}) to the search for critical points of the function $W$
defined in (\ref{eq:1.4}). Observe that the assumption on the
nondegeneracy of the critical point is a rather mild assumption
since it is certainly fulfilled for generic choice of the open
domain $\Omega$.

\medskip

Semilinear equations involving fourth order elliptic operator and
exponential nonlinearity appear naturally in conformal geometry and
in particular in the prescription of the so called $Q$-curvature on
$4$-dimensional Riemannian manifolds \cite{Cha-1}, \cite{Chang-Yang}
\[
Q_g = \frac{1}{12} \, \left( -  \Delta_g S_g + S^2_g - 3 \,
|\mbox{Ric}_g|^2\right)
\]
where $\mbox{Ric}_g$ denotes the Ricci tensor and $S_g$ is the
scalar curvature of the metric $g$. Recall that the $Q$-curvature
changes under a conformal change of metric
\[
g_{w} = e^{2 w} \, g,
\]
according to
\begin{equation}
P_{g} \, w + 2 \, Q_g = 2 \, \tilde{Q}_{g_w} \, e^{ 4 \, w}
\label{eq:1.5}
\end{equation}
where
\begin{equation}
P_{g} : =  \Delta_g^2 +  \delta \, \left( \frac{2}{3}\, S_g \, I - 2
\, \mbox{Ric}_g \right) \, d \label{eq:1.6}
\end{equation}
is the Panietz operator, which is an elliptic $4$-th order partial
differential operator \cite{Chang-Yang} and which transforms
according to
\begin{equation}
e^{4 w} \, P_{e^{2w} g}  =  P_g, \label{eq:1.7}
\end{equation}
under a conformal change of metric $ g_w  : = e^{2w} \, g$. In the
special case where the manifold is the Euclidean space, the Panietz
operator is simply given by
\[
P_{g_{eucl}} = \Delta^{2}
\]
in which case (\ref{eq:1.5}) reduces to
\[
\Delta^{2} \, w =  \tilde{Q} \, e^{4 \, w}
\]
the solutions of which give rise to conformal metric $g_w = e^{2 \,
w} \, g_{eucl}$ whose $Q$-curvature is given by $\tilde Q$. There is
by now an extensive literature about this problem and we refer to
\cite{Chang-Yang} and \cite{Malchiodi-Djadli} for references and
recent developments.

\medskip

When $n=2$, the analogue of the $Q$-curvature reduces is the Gauss
curvature and the corresponding problem has been studied for a long
time. More relevant to the present paper is the study of nontrivial
branches of solutions of
\begin{equation}
\left\{\begin{array}{rclll}-\Delta u &=&
\rho^{2} e^{u} &\mbox{in}& \Omega\\\\u& =&0 & \mbox{on}&
\partial\Omega,
\end{array}
\right. \label{eq:1.8}
\end{equation}
which are defined on some domain of ${\mathbb R}^2$. The study of
this equation goes back $1853$ when Liouville derived a
representation formula for all solutions of (\ref{eq:1.8}) which are
defined in ${\mathbb R}^2$, \cite{Lio}.

\medskip

It turn out that, beside the applications in geometry, elliptic
equations with exponential nonlinearity also arise in the modeling
of many physical phenomenon such as~: thermionic emission,
isothermal gas sphere,  gas combustion, gauge theory
\cite{Tarantello}, \ldots

\medskip

When $\rho$ tends to $0$, the asymptotic behavior nontrivial
branches of solutions of (\ref{eq:1.8}) is well understood thanks
to the pioneer work of Suzuki \cite{Suz} which characterizes the
possible limit of nontrivial branches of solutions of
(\ref{eq:1.8}). The existence of nontrivial branches of solutions
was first proven by Weston \cite{Wes} and then a general result
has been obtained by Baraket and Pacard \cite{Bar-Pac}. More
recently these results were extended, with applications to the
Chern-Simons vortex theory in mind, by Esposito, Grossi and
Pistoia \cite{Esp-Gro-Pis} and Del Pino, Kowalczych and Musso
\cite{del-Mus} to handle equations of the form
\[
-\Delta u  = \rho^2 \, V \, e^u
\]
where $V$ is a non constant (positive) function. We give in Section
9 some results concerning the fourth order analogue of this
equation. Let us also mention that the construction of nontrivial
branches of solutions of semilinear equations with exponential
nonlinearities has allowed Wente to provide counterexamples to a
conjecture of Hopf \cite{Wen} concerning the existence of compact
(immersed) constant mean curvature surfaces in Euclidean space.

\medskip

We now describe the plan of the paper~: In Section 2 we discuss
rotationally symmetric solutions of (\ref{eq:1.1}). In Section 3 we
study the linearized operator about the radially symmetric solution
defined in the previous section. In Section 4, we discuss the
analysis of the bi-Laplace operator in weighted spaces. Both section
strongly use the b-operator which has been developed by Melrose
\cite{Mel} in the context of weighted Sobolev spaces and by Mazzeo
\cite{Maz} in the context of weighted H\"older spaces (see also
\cite{Pac-Riv}).

\medskip

A first nonlinear problem is studied in Section 6  where the
existence of an infinite dimensional family of solutions of
(\ref{eq:1.1}) which are defined on a large ball and which are close
to the rotationally symmetric solution is proven. In Section 7, we
prove the existence of an infinite dimensional family of solutions
of (\ref{eq:1.1}) which are defined on $\Omega$ with small ball
removed. Finally, in Section 8, we show how elements of these
infinite dimensional families can be connected to produce solutions
of (\ref{eq:1.1}) described in Theorem~\ref{th:1.1}. This last
section borrow ideas from applied mathematics were domain
decomposition methods are of common use. In Section 9 is devoted to
some comments. In Section 10, we explain how the results of the
previous analysis can be extended to handle equations of the form
$\Delta^2 u = \rho^4 \, V \, e^u$.

\section{Rotationally symmetric solutions}

We first describe the rotationally symmetric solutions of
\begin{equation}
\Delta^{2}u - \rho^{4}e^{u} = 0 ,\label{eq:2.1}
\end{equation}
which will play a central r\^ole in our analysis. By the
classification given by \cite{Lin}, for $\varepsilon>0$, we define
\[
u_{\varepsilon} (x) : =  4 \, \log(1 + \varepsilon^{2}) - 4 \,
\log(\varepsilon^{2} + |x|^2).
\]
which is clearly a solution of (\ref{eq:2.1}) when
\begin{equation}
\rho^{4} = \frac{{384} \,
{\varepsilon^{4}}}{(1+\varepsilon^{2})^{4}}. \label{eq:2.2}
\end{equation}

Let us notice that equation (\ref{eq:2.1}) is invariant under some
dilation in the following sense~: If $u$ is a solution of
(\ref{eq:2.1}) and if $\tau
> 0$, then $u(\tau \, \cdot ) + 4 \, \log \, \tau$ is also a
solution of (\ref{eq:2.1}). With this observation in mind, we
define, for all $\tau > 0$
\begin{equation}
u_{\varepsilon,\tau} (x) : = 4 \, \log \, (1 + \varepsilon^{2})+ 4
\, \log \, \tau - 4 \, \log \, (\varepsilon^{2} + \tau^{2} \,
|x|^2). \label{eq:2.3}
\end{equation}

\section{A linear fourth order elliptic operator on $\mathbb{R}^{4}$}

We define the linear fourth order elliptic operator
\begin{equation}
{\mathbb L} : = \Delta^{2} - \frac{384}{(1 + |x|^{2})^{4}}
\label{3.1}
\end{equation}
which corresponds to the linearization of
(\ref{eq:2.1}) about the solution $u_1 (= u_{\e =1})$ which has been
defined in the previous section.

\medskip

We are interested in the classification of bounded solutions of
${\mathbb L} \, w =0$ in ${\mathbb R}^4$. Some solutions are easy to
find. For example, we can define
\[
\phi_{0}(x) : = r \, \partial_{r}u_{1}(x) + 4 = 4 \, \frac{1-
r^{2}}{1+r^{2}},
\]
where $r=|x|$. Clearly ${\mathbb L} \, \phi_0 =0$ and this reflects
the fact that (\ref{eq:2.1}) is invariant under the group of
dilations $ \tau \longrightarrow u(\tau \, \cdot ) + 4 \, \log \,
\tau $. We also define, for $i=1, \ldots, 4$
\[
\phi_i(x) : = -{\partial_{x_i}} u_{1}(x) =\frac{8 \,
x_i}{1+|x|^{2}},
\]
which are also solutions of ${\mathbb L}  \, \phi_j =0$ since these
solutions correspond to the invariance of the equation under the
group of translations $a \longrightarrow u( \cdot +a)$.

\medskip

The following result classifies all bounded solutions of ${\mathbb
L} \, w =0$ which are defined in ${\mathbb R}^4$.
\begin{lemma}
Any bounded solution of ${\mathbb L} \, w = 0$ defined in ${\mathbb
R}^4$ is a linear combination of $\phi_i$ for $i=0, 1, \ldots ,4$.
\label{le:3.1}
\end{lemma} {\bf Proof : }  We consider on ${\mathbb R}^4$ the Euclidean metric $g_{eucl} =
dx^2$ and the spherical metric
\[
g_{S^4} =  \frac{4}{(1+|x|^2)^2} \, dx^2
\]
induced by
\[
\begin{array}{rclcllll}
\Pi : & {\mathbb R}^4 & \longrightarrow & S^4 \\[3mm]
& x & \longmapsto & \left( \frac{2x}{1+|x|^2} , \frac{1-|x|^2}{1+
|x|^2} \right)
\end{array}
\]
the inverse of  the stereographic projection.

\medskip

According to \cite{Chang-Yang} we have  $P_{g_{S^4}} =
\Delta_{S^4}^2 - 2 \, \Delta_{S^2}$ and $P_{g_{eucl}} = \Delta^2$.
Therefore, we obtain from (\ref{eq:1.7})
\[
\left(\frac{4}{(1+|x|^2)^2}\right)^2 \, \left( \Delta_{S^4}^2 - 2 \,
\Delta_{S^2} \right)  =  \Delta^2
\]
In particular, if $w : {\mathbb R}^4 \longrightarrow {\mathbb R}$ is
a bounded solution of ${\mathbb L} \, w = 0$ then, $w : S^4-\{N\}
\longrightarrow {\mathbb R}$ is a bounded solution of
\begin{equation}
( \Delta^{2}_{S^4} - 2 \, \Delta_{S^4} - 24  )\, w =0 \label{eq:3.2}
\end{equation}
away from the north p\^ole $N \in S^4$. It is easy to check that the
isolated singularity at the north p\^ole is removable (since $w$ is
assumed to be bounded) and hence (\ref{eq:3.2}) holds on all $S^4$.

\medskip

We now perform the eigenfunction decomposition of $w$ in terms of
the eigendata of the Laplacian on $S^4$. We decompose
\[
w = \sum_{\ell \geq 0} w_\ell
\]
where $w_\ell$ belongs to the $\ell$-th eigenspace of
$\Delta_{S^4}$, namely, $w_\ell$ satisfies $\Delta_{S^4} w_\ell = -
\lambda_\ell \, w_\ell$ with \[ \lambda_\ell : = \ell \, (\ell+3).
\]
We get from (\ref{eq:3.2})
\[
\, (\lambda_\ell^{2}+2 \, \lambda_\ell - 24 ) \, w_\ell = 0 .
\]
Hence, $w_\ell =0$ for all $\ell$ except eventually those for which
$\lambda_\ell = 4$. This implies that $w : S^4 \longrightarrow
{\mathbb R}$ is a combination of the eigenfunctions associated to
$\ell=1$ which are given by $\varphi_i (y) = y_i$ for $i=1, \ldots,
5$, where $y =(y_1, \ldots, y_5) \in S^4$. The sphere being
parameterized by the inverse of the stereographic projection we may
write $y = \Pi (x)$. Then, the functions $4 \, \varphi_i$ precisely
correspond to the functions $\phi_i$ for $i=1, \ldots, 4$, while the
function $4 \, \varphi_5$ corresponds to the function $\phi_0$. This
completes the proof of the result. \hfill $\Box$

\medskip

Let $B_r$ denote the ball of radius $r$ centered at the origin in
${\mathbb R}^4$.
\begin{definition}
Given $k \in {\mathbb N}$, $\alpha \in (0,1)$ and $\mu \in {\mathbb
R}$, we define the H\"older weighted space ${\cal
C}^{k,\alpha}_\mu(\mathbb{R}^{4})$ as the space of functions $w \in
{\cal C}^{k,\alpha}_{loc}(\mathbb{R}^{4})$ for which the following
norm
\[
\| w \|_{{\cal C}^{k,\alpha}_\mu(\mathbb{R}^{4})} : = \|
w\|_{{\mathcal C}^{k, \alpha} (\bar B_1)} +
\displaystyle{\sup_{r\geq 1}} \, \left( (1+r^2)^{-\mu/2} \, \|w
(r\,\cdot) \|_{{\cal C}^{k,\alpha}_\mu(\bar B_1-B_{1/2})} \right),
\]
is finite. \label{de:3.1}
\end{definition}
More details about these spaces and their use in nonlinear problem
can be found in \cite{Pac-Riv}. Roughly speaking, functions in
${\cal C}^{k,\alpha}_\mu(\mathbb{R}^{4})$ are bounded by a
constant times $(1+r^2)^{\mu/2}$ and have $\ell$-th partial
derivatives which are bounded by $(1+r^2)^{(\mu -\ell)/2}$, for
$\ell=1, \ldots, k + \alpha$.

\medskip

As a consequence of the result of Lemma~\ref{le:3.1}, we have the~:
\begin{proposition}
Assume that $\mu > 1$ and $\mu \not\in \mathbb{N}$, then
\[
\begin{array}{rclclllll}
L_\mu : & {\cal C}^{4,\alpha}_\mu(\mathbb{R}^{4}) & \longrightarrow
& {\cal C}^{0,\alpha}_{\mu-4}(\mathbb{R}^{4})\\[3mm]
& w & \longmapsto & {\mathbb L} \, w \end{array}
\]
is surjective. \label{pr:3.1}
\end{proposition}
{\bf Proof :} The mapping properties of $L_\mu$ are very sensitive
to the choice of the weight $\mu$. In particular, it is proved in
\cite{Mel} and \cite{Maz}  (see also \cite{Pac-Riv}) that $L_\mu$
has closed range and is Fredhlom provided $\mu$ is not an indicial
root of $\mathbb L$ at infinity. Recall that $\zeta \in {\mathbb R}$
is an indicial root of ${\mathbb L}$ at infinity if there exists a
smooth function $v$ on $S^3$ such that
\[
{\mathbb L} \, (|x|^\zeta \, v ) = {\mathcal O} (|x|^{\zeta -5})
\]
at infinity.  It is easy to check that the indicial roots of
${\mathbb L}$ at infinity are all $\zeta \in {\mathbb Z}$. Indeed,
let $e$ be an eigenfunction of $\Delta_{S^{3}}$ which is associated
to the eigenvalue $\gamma \, (\gamma +2)$, where $\gamma \in
{\mathbb N}$, hence
\[
\Delta_{S^{3}} \, e = - \gamma \, (\gamma +2 ) \, e .
\]
Then
\[
{\mathbb L} \, (|x|^{\zeta} \, e ) = (\zeta - \gamma) (\zeta -
\gamma - 2) \, (\zeta + 2 +\gamma) \, (\zeta  + \gamma ) \,
|x|^{\zeta -4} \, e  + {\mathcal O}(|x|^{\zeta-8}).
\]
Therefore, we find that $-\gamma - 2$, $-\gamma$, $\gamma$ and
$\gamma+2$ are indicial roots of $\mathbb L$ at $0$. Since the
eigenfunctions of the Laplacian on the sphere constitute a Hilbert
basis of $L^{2}(S^{3})$, we have obtained all the indicial roots of
$\mathbb L$ at infinity.

\medskip

If $\mu \notin {\mathbb Z}$, some duality argument (in weighted
Sobolev spaces) shows that the operator $L_\mu$ is surjective if and
only if the operator $L_{-\mu}$ is injective. And, still under this
assumption \[\mbox{dim}  \, \mbox{Ker}  \, L_\mu  = \mbox{dim} \,
\mbox{Coker} \, L_{-\mu} . \] The result of Lemma~\ref{le:3.1}
precisely states that the operator $L_\mu$ is injective when $\mu <
- 1$.  Therefore, we conclude that $L_\mu$ is surjective when $\mu >
1$, $\mu \not\in \mathbb{Z}$. This completes the proof of the
result. \hfill $\Box$

\section{Analysis of the bi-Laplace operator in weighted spaces}

Given $x^1, \ldots, x^m \in \Omega$ we define $X: = (x^1, \ldots,
x^m)$ and
\[
\bar \Omega^*  \, (X): = \bar \Omega - \{x^1, \ldots , x^m\},
\]
and we choose $r_0 > 0$ so that the balls $B_{r_0}(x^i)$ of center
$x^i$ and radius $r_0$ are mutually disjoint and included in
$\Omega$. For all $r \in (0, r_0)$ we define
\[
\bar \Omega_r \, (X) : = \bar \Omega - \cup_{j=1}^{m} B_{r}(x^j)
\]
With these notations, we have the~:
\begin{definition}
Given $k \in \mathbb{R}$, $\alpha \in (0,1)$ and $\nu\in
\mathbb{R}$, we introduce the H\"older weighted space
${\mathcal{C}}^{k,\alpha}_{\nu}(\bar \Omega^* \, (X))$ as the space
of functions $w \in {\mathcal{C}}^{k,\alpha}_{loc}(\bar \Omega^* \,
(X))$ which is endowed with the norm
\[
\| w \|_{{\mathcal{C}}^{ k, \alpha}_{\nu}(\bar \Omega^* \, (X))}
 : = \| w \|_{{\mathcal{C}}^{ k,\alpha}( \bar \Omega_{r_0/2} \, (X))}+ \sum_{j=1}^m \, \sup_{r
 \in (0, r_0/2)}  \, \left( r^{-\nu} \, \| w(x^j + r \, \cdot )
\|_{{\mathcal{C}}^{ k,\alpha} (\bar B_2-B_1)} \right),
\]
is finite. \label{de:4.1}
\end{definition}
Again, these spaces have already been used many times in nonlinear
contexts and we refer to \cite{Pac-Riv} for further details and
references. Functions which belong to ${\mathcal{C}}^{ k,
\alpha}_{\nu}(\bar \Omega^* \, (X))$ are bounded by a constant
times the distance to $X$ to the power $\nu$ and have $\ell$-th
partial derivatives which are bounded by a constant times the
distance to $X$ to the power $\nu -\ell$, for $\ell =1, \ldots, k
+ \alpha$.

\medskip

When $k \geq 2$, we denote by $[ {\cal C}^{k,\alpha}_\nu (\bar
\Omega^* \, (X)) ]_0$ be the subspace of functions $w \in {\cal
C}^{k,\alpha}_\nu (\bar \Omega^* \, (X))$ satisfying $ w = \Delta w
= 0$.

\medskip

We will use the following~:
\begin{proposition}
Assume that $\nu < 0$ and $\nu \not\in \mathbb{Z}$, then \[
\begin{array}{rclcllll}
 {\mathcal L}_\nu & : [ {\cal C}^{4,\alpha}_\nu (\bar \Omega^* \, (X))]_0 &
\longrightarrow  & {\cal C}^{0,\alpha}_{\nu-4} (\bar \Omega^* \, (X))\\[3mm]
& w & \longmapsto & \Delta^2 \, w
\end{array}
\]
is surjective. \label{pr:4.1}
\end{proposition} {\bf Proof~:}
Again this result follows from the theory developed in \cite{Mel}
and \cite{Maz} (see also \cite{Pac-Riv}). The mapping properties of
${\mathcal L}_\nu$ depend on the choice of the weight $\nu$. The
operator ${\mathcal L}_\nu$ has closed range and is Fredhlom
provided $\nu$ is not an indicial root of $\Delta^2$ at the points
$x^j$. Recall that $\zeta \in {\mathbb R}$ is an indicial root of
${\mathcal L}_\nu$ at $x^j$ if there exists a smooth function $v$ on
$S^3$ such that
\[
{\mathbb L} \, (|x-x^j|^\zeta \, v ) = {\mathcal O} (|x-x^j|^{\zeta
-3})
\]
at $x^j$. As in Proposition~\ref{pr:4.1}, it is easy to check that
the indicial roots of ${\mathcal L}$ at $x^j$ are all $\zeta \in
{\mathbb Z}$.

\medskip

If $\nu \notin {\mathbb Z}$, some duality argument (in weighted
Sobolev spaces) shows that the operator ${\mathcal L}_\nu$ is
surjective if and only if the operator ${\mathcal L}_{-\nu}$ is
injective. And, still under this assumption \[ \mbox{dim}  \,
\mbox{Ker}  \, {\mathcal L}_\nu  =  \mbox{dim} \, \mbox{Coker} \,
{\mathcal L}_{-\nu} .
\]

We claim that the operator ${\mathcal L}_\nu$ is injective if $\nu >
0 $. Indeed, isolated singularities of any solution $w \in {\mathcal
C}^{4, \alpha}_\nu  (\bar \Omega^* \, (X))$ of $\Delta^2 w =0$ in
$\Omega^*$ are removable if $\nu > 0$. Therefore, $w$ is a
bi-harmonic function in $\Omega$ with $w = \Delta w=0$ on $\partial
\Omega$. This implies that $w \equiv 0$ and hence ${\mathcal L}_\nu$
is injective when $\nu >0$ as claimed.

\medskip

We then conclude that ${\mathcal L}_\nu$ is surjective when $\nu <
0$, $ \nu \notin {\mathbb Z}$. This completes the proof of the
result. \hfill $\Box$

\medskip

Given $y^1, \ldots, y^m$ close enough to $x^1, \ldots, x^m$, we set
$Y : = (y^1, \ldots, y^m)$ and we define a family of diffeomorphisms
$D (=  D_{X,Y})$
\[
D : \Omega \longrightarrow \Omega
\]
depending smoothly on $y^1, \ldots, y^m$ by
\begin{equation}
D (x) = x + \sum_{j=1}^m \chi_{r_0} (x-x^j) \, (x^j- y^j),
\label{eq:4.1}
\end{equation}
where $\chi_{r_0}$ is a cutoff function identically equal to $1$ in
$B_{r_0/2}$ and identically equal to $0$ outside $B_{r_0}$. In
particular, $D(y^j) = x^j$ for each $j$.

\medskip

The equation $\Delta^2 \, \tilde w = \tilde f$ where  $\tilde f \in
{\mathcal C}^{0, \alpha}_{\nu-4} (\bar \Omega \, (Y))$ can be solved
by writing $\tilde w  = w \circ D$ and $\tilde f  = f\circ D$ so
that $w$ is a solution of the problem
\begin{equation}
\Delta^2 \, w  + \left( \Delta^2 (w \circ D)  - (\Delta^2 w )\circ
D \right) \circ D^{-1} = f \label{eq:4.2}
\end{equation}
where this time  $f  \in {\mathcal C}^{0, \alpha}_\nu (\bar \Omega
 \, (X))$. It should be clear that
\begin{equation}
\| \left( \Delta^2 (w \circ D)  - (\Delta^2 w )\circ D \right) \circ
D^{-1} \|_{{\mathcal C}^{0, \alpha}_{\nu -4} (\bar \Omega^* \, (X))}
\leq c \, \| Y-X\|  \, \| w\|_{{\mathcal C}^{4, \alpha}_\nu (\bar
\Omega^* \, (X))} \label{eq:4.3}
\end{equation}
provided $\| Y - X \| \leq r_0/4$.

\medskip

We fix $\nu <0$, $\nu \notin {\mathbb Z}$ and use the result of
Proposition~\ref{pr:4.1} to choose a right inverse ${\mathcal
G}_{\nu, X}$  for ${\mathcal L}_\nu : [{\mathcal C}^{4, \alpha}_\nu
(\bar \Omega^* \, (X))]_0 \longrightarrow {\mathcal C}^{0,
\alpha}_{\nu -4} (\bar \Omega^* \, (X))$. The estimate
(\ref{eq:4.3}) together with a perturbation argument, shows that
(\ref{eq:4.2}) is solvable provided $Y$ is close enough to $X$. This
provides a right inverse ${\mathcal G}_{\nu, Y}$ which depends
continuously (and in fact smoothly) on the points $y^1, \ldots, y^m$
in the sense that
\[
f \quad \longmapsto \quad {\mathcal G}_{\nu, Y} (f \circ D_{X,Y} )
\circ (D_{X,Y})^{-1}
\]
depends smoothly on $Y$.

\section{Bi-harmonic extensions}

Given $\varphi \in {\cal C}^{4,\alpha}(S^3)$ and  $\psi \in {\cal
C}^{2,\alpha}(S^3)$ we define $H^i (= H^i (\varphi, \psi \, ; \cdot
) )$ to be the solution of
\begin{equation}\left\{\begin{array}{rclll}\Delta^{2} \, H^i&=& 0
&\mbox{in}& B_1\\[3mm]
H^i & =& \varphi  &\mbox{on}& \partial B_1\\[3mm]
\Delta \, H^i & =& \psi &\mbox{on}& \partial B_1,
\end{array}
\right. \label{eq:5.1}
\end{equation}
where, as already mentioned, $B_1$ denotes the unit ball in
${\mathbb R}^4$.

\medskip

We set $B_1^* = B_1 - \{0\}$. As in the previous section, we
define~:
\begin{definition}
Given $k\in {\mathbb N}$, $\alpha \in (0,1)$ and $\mu \in {\mathbb
R}$, we introduce the H\"older weighted spaces ${\cal
C}^{k,\alpha}_\mu(\bar B_1^*)$ as the space of function in ${\cal
C}^{k,\alpha}_{loc}(\bar B_1^*)$ for which the following norm
\[
\|u \|_{{\cal C}^{k,\alpha}_\mu(\bar B_1^*)} =
\displaystyle{\sup_{r\leq 1/2}} \, \left( r^{-\mu} \, \| u(r \,
\cdot) \|_{{\cal C}^{k,\alpha}(\bar B_2-B_1)} \right),
\]
is finite. \label{de:5.1}
\end{definition} This corresponds to the
space and norm already defined in the previous section when $\Omega
= B_1$, $m=1$ and $x^1=0$.

\medskip

Let $e_1, \ldots, e_4$ be the coordinate functions on $S^3$. We
prove the~:
\begin{lemma}
Assume that
\begin{equation}
\int_{S^{3}} (8 \, \varphi  - \psi ) \,  dv_{S^3} =0 \qquad
\mbox{and also that} \qquad \int_{S^{3}} (12 \, \varphi - \psi ) \,
e_\ell \, dv_{S^3}=0 \label{eq:5.2}
\end{equation}
for $\ell = 1, \ldots , 4.$ Then there exists $c> 0$ such that
\[
\| H^i(\varphi, \psi \, ; \cdot ) \|_{{\mathcal{C}}_{2}^{ 4,\alpha}(
\bar B_{1}^*)}\leq c \, ( \| {\varphi} \|_{{\mathcal{C}}^{
4,\alpha}( S^3)} + \| {\psi} \|_{{\mathcal{C}}^{ 2,\alpha}(S^3)} ).
\]
\label{le:5.1}
\end{lemma}
{\bf Proof~:} There are many ways to proof this result. Here is a
simple one which has the advantage to be quite flexible. We consider
the eigenfunction decomposition of $\varphi$ and $\psi$ in terms of
the eigenfunctions of $\Delta_{S^3}$.
\begin{equation}
\varphi = \sum_{\ell \geq 0}  \varphi_\ell  \qquad\mbox{and}\qquad
\psi =\sum_{\ell \geq 0}  \psi_\ell , \label{eq:5.3}
\end{equation}
where, for each $\ell \geq 0$, the functions $\varphi_\ell$ and
$\psi_\ell$ belong to the $\ell$-th eigenspace of $\Delta_{S^3}$,
namely
\[
\Delta_{S^3} \varphi_\ell = - \ell \, (2+\ell) \, \varphi_\ell
\qquad \mbox{and} \qquad \Delta_{S^3} \psi_\ell =  - \ell \,
(2+\ell) \, \psi_\ell .
\]
Then the function $H^i$ can be explicitly written as
\begin{equation}
H^i  = \sum_{\ell \geq 0}   \, r^\ell \, \left( \varphi_\ell -
\frac{1}{4(\ell+2)} \, \psi_\ell \right) + \sum_{\ell \geq 0}
\frac{1}{4(\ell+2)} r^{2+\ell} \, \psi_\ell. \label{eq:5.4}
\end{equation}
Observe that, under the hypothesis, the coefficients of $r^0$ and
$r^1$ vanish and hence, at least formally, the expansion of $H$ only
involves powers of $r$ which are greater than or equal to $2$.

\medskip

We claim that
\[
\| \varphi_\ell \|_{L^\infty} \leq c_\ell \, \| \varphi \|_{L^2}
\qquad \qquad \| \psi_\ell \|_{L^\infty} \leq c_\ell \, \| \psi
\|_{L^2}
\]
where the constant $c_\ell$ depends polynomially on $\ell$. For
example, we can write $ \varphi_\ell =  a_\ell \, e_\ell$ where
$a_\ell \in {\mathbb R}$ and $e_\ell$ is an eigenfunction  of
$\Delta_{S^3}$ which is normalized to have $L^2$ norm equal to
$1$. Then \[ |a_\ell |  =
 \left| \int_{S^3} \, \varphi_\ell \, e_\ell \, dv_{S^3} \right| \leq
 c \, \| \varphi \|_{L^2} \leq c \, \| \varphi\|_{L^\infty}
\]
Next, $e_\ell$ solves $\Delta_{S^3} e_\ell = -\ell \, (2+\ell) \,
e_\ell $, we can use elliptic regularity theory to show that the
$L^\infty (S^3)$ norm of $e_\ell$ depends polynomially on $\ell$.
The claim then follows at once.

\medskip

This immediately yields the estimate
\[
\sup_{r\leq 1/2}  \, \left( r^{-2} \, |H^i| + |\Delta \, H^i|
\right) \leq c \, (\| \varphi\|_{L^\infty} + \|\psi \|_{L^\infty})
\]
This estimate, together with the maximum principle and standard
elliptic estimates yields
\[
\sup_{r\leq 1} r^{-2} \, |H^i|\leq c \, (\| \varphi\|_{L^\infty} +
\|\psi \|_{L^\infty})
\]
The estimate for the derivatives of $H^i$ now follows at once from
Schauder's estimates. \hfill $\Box$

\medskip

Given $\varphi \in {\cal C}^{4,\alpha}(S^3)$ and $\psi \in {\cal
C}^{2,\alpha}(S^3)$ we define (when it exists !) $H^e( =H^e
(\varphi, \psi \, ; \cdot))$ to be the solution of
\begin{equation}
\left\{\begin{array}{rclll}\Delta^{2} \, H^e &=& 0
& \mbox{in} & {\mathbb R}^4 - B_1\\[3mm]
H^e & = & \varphi & \mbox{on}& \partial B_1\\[3mm]
\Delta H^e & =& \psi &\mbox{on}& \partial B_1,
\end{array}
\right. \label{eq:5.5}
\end{equation}
which decays at infinity.
\begin{definition}
Given $k\in {\mathbb N}$, $\alpha \in (0,1)$ and $\nu \in {\mathbb
R}$, we define the space ${\cal C}^{k,\alpha}_\nu ({\mathbb R}^4
-B_1)$ as the space of functions  $w \in {\cal
C}^{k,\alpha}_{loc}({\mathbb R}^4 -B_1)$ for which the following
norm
\[
\| w \|_{{\cal C}^{k,\alpha}_\nu({\mathbb R}^4 -B_1)} =
\displaystyle{\sup_{r\geq 1}} \, \left(  r^{-\nu} \, \|w(r \,
\cdot)\|_{{\cal C}^{k,\alpha}_\nu(\bar B_2-B_1)} \right),
\]
is finite. \label{de:5.2}
\end{definition}
\medskip

We prove the~:
\begin{lemma}
Assume that
\begin{equation}
\int_{S^{3}} \psi  \,  dv_{S^3} = 0.\label{eq:5.6}
\end{equation}
Then there exists $c > 0$ such that \[ \| H^e (\varphi, \psi \, ;
\cdot) \|_{{\mathcal{C}}_{-1}^{ 4,\alpha} ({\mathbb R}^4 -B_1)}\leq
c \, (\|{\varphi}\|_{{\mathcal{C}}^{ 4,\alpha}( S^3)}+
\|{\psi}\|_{{\mathcal{C}}^{ 2,\alpha}(  S^3)} ).
\]
\label{le:5.2}
\end{lemma}
{\bf Proof~:} We use the notations of the previous Lemma. Now, the
function $H^e$ can be explicitly written as
\begin{equation}
H^e  = r^{-2} \,\varphi_0 + \sum_{\ell \geq 1}   \, r^{-2-\ell} \,
\left( \varphi_\ell + \frac{1}{4 \, \ell} \, \psi_\ell \right) -
\sum_{\ell \geq 1} \frac{1}{4 \, \ell} \, r^{-\ell} \, \psi_\ell.
\label{eq:5.7} \end{equation} Observe that,  (\ref{eq:5.6}) implies
that the coefficients of $r^0$, vanishes and hence the expansion of
$H^e$ only involves powers of $r$ which are lower than or equal to
$-1$. The proof is now identical to the proof of Lemma~\ref{le:5.1}
and left to the reader. \hfill $\Box$

\medskip

Under the hypothesis of the Lemma \ref{le:5.1}, there is uniqueness
of the bi-harmonic extension of the boundary data which decays at
infinity.

\medskip

If $F\subset L^2(S^3)$ is a space of functions defined on $S^3$, we
define the space $F^\perp$ to be the subspace of functions of $F$
which are $L^2(S^3)$-orthogonal to the functions $1, e_1, \ldots,
e_4$. We will need the~:
\begin{lemma}
The mapping
\[
\begin{array}{lclclllll}
\mathcal{P}: & \mathcal{C}^{4,\alpha}(S^3)^\perp \times
\mathcal{C}^{2,\alpha}(S^3)^\perp  & \longrightarrow &
\mathcal{C}^{3,\alpha}(S^3)^\perp \times \mathcal{C}^{1,\alpha}(S^3)^\perp  \\[3mm]
 &(\,\varphi ,\psi\,) & \longmapsto & (\partial_{r}H^i - \del_r H^e  \, , \, \partial_{r} \, \Delta H^i
 - \del_r \, \Delta H^e \,)
\end{array}
\]
where $H^i = H^i ( \varphi , \psi \, ; \cdot )$ and $H^e = H^e(
\varphi , \psi \, ; \cdot )$, is an isomorphism. \label{le:5.3}
\end{lemma}
{\bf Proof~:} Granted the explicit formula given in the previous two
Lemmas, we have
\begin{equation}
{\mathcal P}  (\varphi, \psi) = \left( \sum_{\ell \geq 2} (\ell+1)
\, \left( 2 \, \varphi_\ell + \frac{1}{\ell(\ell+2)} \, \psi_\ell
\right) , \sum_{\ell \geq 2} 2 \, (\ell+1) \, \psi_\ell \right).
\label{eq:5.8}
\end{equation}
 We denote by $W^{k, 2} (S^3)$ Sobolev
space of functions whose weak partial derivatives, up to order $k$
are in $L^2 (S^3)$. The norm in $W^{k, 2} (S^3)$ can be chosen to be
\[
\| \varphi  \|_{W^{k, 2} (S^3)} : = \left( \sum_{\ell \geq 0}
(1+\ell)^{2k} \, \|\varphi_\ell \|_{L^2 (S^3)}^2 \right)^{1/2}
\]
when the function $\varphi$ is decomposed over eigenspaces of
$\Delta_{S^3}$
\[
\varphi = \sum_{\ell \geq 0} \, \varphi_\ell
\]
where $\Delta_{S^3} \, \varphi_\ell =  - \ell \, ( \ell+2) \,
\varphi_\ell$. It follows at once that
\[
{\mathcal P} : W^{k+3,2}(S^3)^\perp \times W^{k+1,2}(S^3)^\perp
\longrightarrow W^{k+2,2}(S^3)^\perp \times W^{k,2} (S^3)^\perp
\]
is invertible. Elliptic regularity theory then implies that the
corresponding map is also invertible when defined between the
corresponding H\"older spaces. \hfill $\Box$

\section{The first nonlinear Dirichlet problem}

For all $\e , \tau > 0 $, we set
\[
R_\varepsilon : = \tau \;/\sqrt{\varepsilon}.
\]
Given $\varphi \in {\mathcal C}^{4, \alpha}(S^3)$ and $\psi \in
{\mathcal C}^{2, \alpha}(S^3)$ satisfying (\ref{eq:5.2}), we define
\[
{\bf u}: = u_1 + H^i ({ \varphi}, { \psi} \, ; \, ( \cdot /
R_\varepsilon)).
\]
We would like to find a solution $u$ of
\begin{equation}
\Delta^2 \, u - 24 \, e^u = 0  \label{eq:6.1}
\end{equation}
which is defined in $B_{R_\e}$ and which is a perturbation of ${\bf
u}$. Writing $u = {\bf u} + v$, this amounts to solve the equation
\begin{equation}
{\mathbb L} \, v = \frac{384}{(1 + r^2)^4}( e^{H^{i} ({ \varphi}, {
\psi} \, ; \, (\cdot / R_\varepsilon )) + v} - 1 - v),
\label{eq:6.2}
\end{equation}
since $H^i$ is bi-harmonic.

\medskip

We will need the following~:
\begin{definition}
Given $ \bar r \geq 1$, $k \in {\mathbb N}$, $\alpha \in (0,1)$ and
$\mu \in {\mathbb R}$, the weighted space ${\mathcal C}^{k,
\alpha}_{\mu} (B_{\bar r})$ is defined to be the space of functions
$w \in {{\mathcal C}}^{k, \alpha} (B_{\bar r})$ endowed with the
norm
\[
\| w \|_{{{\mathcal C}}^{k, \alpha}_{\mu} (\bar B_{\bar r})}  : = \|
w \|_{{{\mathcal C}}^{k, \alpha} (B_{1})} + \sup_{1 \leq r \leq \bar
r}  \, \left( r^{-\mu} \, \| w (r \, \cdot) \|_{{{\mathcal C} }^{k,
\alpha} ( \bar B_{1} - B_{1/2})} \right).
\]
\label{de:6.1}
\end{definition}

For all $\sigma \geq 1$, we denote by \[{\mathcal E}_{\sigma} :
{\mathcal C}^{0, \alpha}_\mu (\bar B_{\sigma}) \longrightarrow
{\mathcal C}^{0, \alpha}_\mu ({\mathbb R}^4) \] the extension
operator defined by
\[
{\mathcal E}_{\sigma} \, (f) (x) = \chi \left( \frac{|x|}{\sigma}
\right) \, f \left( \sigma \, \frac{x}{|x|} \right)
\]
where $t \longmapsto \chi (t)$ is a smooth nonnegative cutoff
function identically equal to $1$ for $t \geq 2$ and identically
equal to $0$ for $t \leq 1$. It is easy to check that there exists a
constant $c = c (\mu)
>0$, independent of $\sigma \geq 1$, such that
\begin{equation}
\| {\mathcal E}_{\sigma} (w ) \|_{{\mathcal C}^{0, \alpha}_{\mu}
({\mathbb R}^4)} \leq \, c \, \| w \|_{{\mathcal C}^{0,
\alpha}_{\mu} (\bar B_{\sigma})} . \label{eq:6.3}
\end{equation}

We fix
\[
\mu \in (1,2)
\]
and denote by ${\mathcal G}_\mu$ a right inverse provided by
Proposition~\ref{pr:3.1}. To find {\underline a}  solution of
(\ref{eq:6.2}), it is enough to find $v \in {\mathcal C}^{4,
\alpha}_\mu ({\mathbb R}^4)$ solution of
\begin{equation}
v = N( \e, \tau, \varphi, \psi \, ; \, v) \label{eq:6.4}
\end{equation}
where we have defined
\[
 N (\e, \tau, \varphi, \psi \, ; \, v) : = {\mathcal G}_\mu \circ {\mathcal E}_{R_\varepsilon} \left(
\frac{384}{(1 + |\cdot |^2)^4} \, \left( e^{H^{i} ({ \varphi}, {
\psi} \, ; \, ( \cdot / R_\varepsilon)) + v} - 1 - v \right) \right)
\]

Given $\kappa >1$ (whose value will be fixed later on), we now
further assume that the functions $\varphi \in {\mathcal C}^{4,
\alpha} (S^3)$, $\psi \in {\mathcal C}^{2, \alpha} (S^3)$ and the
constant $\tau >0$ satisfy
\begin{equation}
|\log (\tau / \tau_*)|\leq \kappa \, \e  \log 1/\e , \qquad \qquad
\|\varphi \|_{{\cal C}^{4, \alpha} (S^3)}\leq \kappa \, \varepsilon
\qquad \mbox{and} \qquad \| \psi \|_{{\cal C}^{2, \alpha} (S^3)}\leq
\kappa \, \varepsilon , \label{eq:6.5}
\end{equation}
where $\tau_* >0$ is fixed.

\medskip

We have the following technical~:
\begin{lemma}
Given $\kappa > 0$. There exist $\e_\kappa >0$, $c_\kappa
>0$ and $\bar c_\kappa >0$ such that, for all $\e \in (0,
\e_\kappa)$
\begin{equation}
\| N (\e, \tau, \varphi, \psi \, ; \, 0) \|_{{\mathcal C}^{4,
\alpha}_\mu ({\mathbb R}^4)}\leq c_\kappa \, \varepsilon^2.
\label{eq:6.6}
\end{equation}
Moreover,
\begin{equation}
\| N (\e, \tau, \varphi, \psi \, ; \, v_2) - N (\e, \tau, \varphi,
\psi \, ; \, v_1) \|_{{\mathcal C}^{4, \alpha}_\mu ({\mathbb
R}^4)}\leq \bar c_\kappa \, \varepsilon^2 \, \| v_2 - v_1
\|_{{\mathcal C}^{4, \alpha}_\mu ({\mathbb R}^4)} \label{eq:6.7}
\end{equation}
and
\begin{equation}
\begin{array}{rllllll}
\| N (\e, \tau, \varphi_2 ,  \psi_2 \, ; \, v) - N (\e, \tau,
\varphi_1, \psi_1 \, ; \,  v) \|_{{\mathcal C}^{4, \alpha}_\mu
({\mathbb R}^4)} \qquad \qquad \qquad \qquad \\[3mm]
\leq \bar c_\kappa \,
\varepsilon^3 \, \left( \| \varphi_2 - \varphi_1 \|_{{\mathcal
C}^{4, \alpha}(S^3)} + \| \psi_2 - \psi_1 \|_{{\mathcal C}^{2,
\alpha}(S^3)} \right)
\end{array}
\label{eq:6.8}
\end{equation}
provided $\tilde v= v, v_1, v_2 \in {\mathcal C}^{4, \alpha}_\mu
({\mathbb R}^4)$, $\tilde \varphi = \varphi, \varphi_1, \varphi_2
\in {\mathcal C}^{4, \alpha} (S^3)$, $\tilde \psi = \psi, \psi_1,
\psi_2 \in {\mathcal C}^{4, \alpha} (S^3)$ satisfy
\[
\| \tilde v \|_{{\mathcal C}^{4, \alpha}_\mu ({\mathbb R}^4)} \leq 2
\, c_\kappa \, \e^2 , \qquad \|\tilde \varphi \|_{{\cal C}^{4,
\alpha} (S^3)}\leq \kappa \, \varepsilon , \qquad  \| \tilde \psi
\|_{{\cal C}^{2, \alpha} (S^3)}\leq \kappa \, \varepsilon ,
\]
and $|\log (\tau / \tau_*)|\leq \kappa \, \e  \log 1/\e$.
\label{le:6.1}
\end{lemma}
{\bf Proof~:} The proof of these estimates follows from the result
of Lemma~\ref{le:5.1} together with the assumption on the norms of
$\varphi$ and $\psi$. Let $c^{(i)}_\kappa$ denote constants which
only depend on $\kappa$ (provided $\e$ is chosen small enough).

\medskip

It follows from Lemma~\ref{le:5.1} that
\[
\| H^i (\varphi, \psi \, ; \,  \cdot / R_\e) \|_{{\mathcal C}^{4,
\alpha} _{2} ( \bar B_{R_\e})} \leq  c\, \, R_\e^{-2} \, (\| \varphi
\|_{{\mathcal C}^{4, \alpha}(S^3)} + \| \psi \|_{{\mathcal C}^{2,
\alpha}(S^3)}) \leq c^{(1)}_\kappa \, \e^2
\]
Therefore, we get
\[
\left\| (1+ |\cdot|^2)^{-4} \, \left( e^{H^i(\varphi, \psi \, ; \,
\cdot / R_\e)} -1 \right) \right\|_{{\mathcal C}^{0, \alpha}_{\mu-4}
(\bar B_{R_\e})} \leq c^{(2)}_\kappa  \, \e^2
\]
Making use of Proposition~\ref{pr:3.1} together with (\ref{eq:6.3})
we conclude that
\[
\| N (\e, \tau, \varphi, \psi \, ; \,  0) \|_{{\mathcal C}^{4,
\alpha}_\mu ({\mathbb R}^4)}\leq c_\kappa \, \varepsilon^2.
\]

To derive the second estimate, we use the fact that
\[ \left\| (1+ |\cdot|^2)^{-4}
\, e^{H^i({\varphi, \psi}; \cdot / R_\e) } \, \left( e^{v_2} -
e^{v_1} - v_2 + v_1 \right) \right\|_{{\mathcal C}^{0,
\alpha}_{\mu-4} (\bar B_{R_\e})} \leq c^{(3)}_\kappa  \, \e^2 \, \|
v_2 - v_1 \|_{{\mathcal C}^{4, \alpha}_\mu ({\mathbb R}^4)}
\]
and
\[
\left\| (1+ |\cdot|^2)^{-4} \,  \left( e^{H^i(\varphi, \psi \, ; \,
\cdot / R_\e)  }  -1 \right) \, \left(v_2  - v_1 \right)
\right\|_{{\mathcal C}^{0, \alpha}_{\mu-4} (\bar B_{R_\e})} \leq
c^{(4)}_\kappa  \, \e^2\, \| v_2 - v_1 \|_{{\mathcal C}^{4,
\alpha}_\mu ({\mathbb R}^4)} ,
\]
provided $v_1, v_2 \in {\mathcal C}^{4, \alpha}_\mu ({\mathbb R}^4)$
satisfy  $\| v_i \|_{{\mathcal C}^{4, \alpha}_\mu ({\mathbb R}^4)}
\leq 2 \, c_\kappa \, \e^2$.

\medskip

Finally, in order to derive the third estimate, we use
\[
\begin{array}{rllll}
\left\| (1+ |\cdot|^2)^{-4} \,  \left( e^{H^i (\varphi_2, \psi_2 \,
; \, \cdot / R_\e) }  - e^{H^i(\varphi_1 , \psi_1 \, ; \,  \cdot /
R_\e)} \right) \, e^v \right\|_{{\mathcal C}^{0, \alpha}_{\mu-4}
(\bar B_{R_\e})} \qquad \qquad \qquad \\[3mm]
\leq c^{(5)}_\kappa  \, \e^2 \, \| H^i (\varphi_2 - \varphi_1 ,
\psi_2 - \psi_1 \, ; \,  \cdot / R_\e ) \|_{{\mathcal C}^{4,
\alpha}_2 (\bar B_{R_\e})} , \end{array}
\] provided $v \in {\mathcal C}^{4,
\alpha}_\mu ({\mathbb R}^4)$ satisfies $\| v \|_{{\mathcal C}^{4,
\alpha}_\mu ({\mathbb R}^4)} \leq 2 \, c_\kappa \, \e^2$. The
second and third estimates again follows from
Proposition~\ref{pr:3.1} and (\ref{eq:6.3}).  \hfill $\Box$

\medskip

Reducing $\e_\kappa >0$ if necessary, we can assume that,
\begin{equation}
\bar c_\kappa \, \e^{2}  \leq \frac{1}{2} \label{eq:6.9}
\end{equation} for all $\e \in (0, \e_\kappa )$. Then,
(\ref{eq:6.6}) and (\ref{eq:6.7}) in Lemma~\ref{le:6.1} are enough
to show that
\[
v  \longmapsto N (\e, \tau, \varphi, \psi \, ; \, v)
\]
is a contraction from
\[
\{ v \in {\mathcal C}^{4, \alpha}_\mu ({\mathbb R}^4)  \qquad :
\qquad \|v \|_{{\mathcal C}^{4, \alpha}_\mu ({\mathbb R}^4)} \leq 2
\, c_\kappa \, \e^{2} \}
\]
into itself and hence has a unique fixed point $ v (\e, \tau,
\varphi, \psi \, ; \, \cdot )$ in this set. This fixed point is
{\underline a} solution of (\ref{eq:6.2}) in $B_{R_\e}$.

\medskip

Reducing $\e_\kappa$ if this is necessary, it follows from
(\ref{eq:6.7}) and (\ref{eq:6.8}) in Lemma~\ref{le:6.1} that
\begin{equation}
\| v (\e, \tau, \varphi_2, \psi_2 \, ; \, \cdot ) - v (\e, \tau,
\varphi_1, \psi_1 \, ; \, \cdot) \| \leq 2 \, \bar c_\kappa \,
\varepsilon^3 \, \left( \| \varphi_2 - \varphi_1 \|_{{\mathcal
C}^{4, \alpha}(S^3)} + \| \psi_2 - \psi_1 \|_{{\mathcal C}^{2,
\alpha}(S^3)} \right). \label{eq:6.10}
\end{equation}

We summarize this in the~:
\begin{proposition}
Given $\kappa >1$, there exist $\e_\kappa >0$ and $c_\kappa >0$
(only depending on $\kappa$) such that given $\varphi \in
{\mathcal C}^{4, \alpha} (S^3)$, $\psi \in {\mathcal C}^{2,
\alpha} (S^3)$ and $\tau >0$ satisfying (\ref{eq:5.2})  and
\begin{equation}
|\log (\tau / \tau_*)|\leq \kappa \, \e  \log 1/\e , \qquad \qquad
\| \varphi \|_{{\cal C}^{4, \alpha} (S^3)}\leq \kappa \, \varepsilon
\qquad \mbox{and} \qquad \| \psi \|_{{\cal C}^{2, \alpha} (S^3)}\leq
\kappa \, \varepsilon , \label{eq:6.100}
\end{equation}
the function
\[
u (\e, \tau, \varphi, \psi \, ; \, \cdot ) : = u_1 + H^i ({
\varphi}, { \psi} \, ; \,  \cdot / R_\varepsilon) + v(\e, \tau,
\varphi, \psi \, ; \, \cdot ),
\]
solves (\ref{eq:6.1}) in $B_{R_\e}$. In addition
\begin{equation}
\| v (\e, \tau, \varphi, \psi \, ; \, \cdot ) \|_{{\cal
C}^{4,\alpha}_{\mu}({\mathbb R}^4)} \leq 2 \, c_\kappa \, \e^{2}.
\label{eq:6.11} \end{equation} \label{pr:6.1}
\end{proposition}

Observe that the function $ v (\e, \tau, \varphi, \psi \, ; \, \cdot
)$ being obtained as a fixed point for contraction mapping, it
depends continuously on the parameter $\tau$.

\section{The second nonlinear Dirichlet problem}

For all $\e \in (0, r_0^2)$, we set
\[
r_\e = \sqrt \varepsilon .
\]
Recall that $G ( x , \cdot)$ denotes the unique solution of
\[
\Delta^2 \, G ( x , \cdot)  = 64 \, \pi^2 \, \delta_{x}
\]
in $\Omega$, with $G ( x, \cdot) = \Delta \, G (x, \cdot) =0$ on
$\del \Omega$. In addition, the following decomposition holds
\[
G (x, y) = - 8 \, \log |x- y| + R (x, y)
\]
where $ y \longmapsto R(x, y)$ is a smooth function.

\medskip

Given $x^1, \ldots, x^m \in \Omega$. The data we will need are the
following~:

\begin{itemize}

\item[(i)] Points $Y := (y^1, \ldots, y^m )\in \Omega^m$ close enough
to $X : = (x^1, \ldots, x^m)$.

\item[(ii)] Parameters $ \Lambda : = (\lambda^1, \ldots, \lambda^m )\in
{\mathbb R}^m$ close to $0$.

\item[(iii)] Boundary data $\Phi : = (\varphi^1,
\ldots, \varphi^m ) \in ({\mathcal C}^{4, \alpha} (S^3))^m$ and
$\Psi  : = (\psi^1, \ldots, \psi^m ) \in ({\mathcal C}^{2, \alpha}
(S^3))^m$ each of which satisfies (\ref{eq:5.6}).
\end{itemize}

With all these data, we define
\begin{equation} {\bf \tilde u}  : =
\sum_{j=1}^m (1 + \lambda^j ) \, G ( y^j , \, \cdot \, ) +
\sum_{j=1}^m \chi_{r_0} (\cdot - y^j) \, H^{e} ( \varphi^j, \psi^j
\, ; \, (\cdot - y^j) / r_\e ) \label{eq:7.1}
\end{equation}
where  $\chi_{r_0}$ is a cutoff function identically equal to $1$ in
$B_{r_0/2}$ and identically equal to $0$ outside $B_{r_0}$.

\medskip

We define $\rho >0$ by
\[
\rho^4 = \frac{384 \, \e^4}{(1+\e^2)^4} .
\]
We would like to find a solution of the equation
\begin{equation}
\Delta^2 \, u - \rho^{4} \, e^u = 0 ,\label{eq:7.2}
\end{equation}
which is defined in $\bar \Omega_{r_\e} (Y)$ and which is a
perturbation of ${\bf \tilde u}$. Writing $u = {\bf \tilde u} +
\tilde v$, this amounts to solve
\begin{equation}
\Delta^2  \, \tilde v = \rho^{4} \, e^{{\bf \tilde u}  + \tilde{v}}
- \Delta^{2} \, {\bf \tilde u} . \label{eq:7.3}
\end{equation}

We need to define an auxiliary weighed space~:
\begin{definition}
Given $\bar r \in (0, r_0/2)$, $k \in \mathbb{R}$, $\alpha \in
(0,1)$ and $\nu\in \mathbb{R}$, we define the H\"older weighted
space ${\mathcal{C}}^{k,\alpha}_{\nu}(\bar \Omega_{\bar r} \, (X))$
as the space of functions $w \in {\mathcal{C}}^{k,\alpha}(\bar
\Omega_{\bar r} \, (X))$ which is endowed with the norm
\[
\| w \|_{{\mathcal{C}}^{ k, \alpha}_{\nu}(\bar \Omega_{\bar r} \,
(X))} : = \| w \|_{{\mathcal{C}}^{ k,\alpha}( \bar \Omega_{r_0/2}
 \, (X) )}+ \sum_{j=1}^m \, \sup_{r \in [\bar r, r_0/2)}  \, \left(
r^{-\nu} \, \| w(x^j + r \, \cdot ) \|_{{\mathcal{C}}^{ k,\alpha}
(\bar B_2-B_1)} \right).
\]
\label{de:7.1}
\end{definition}

For all $\sigma \in (0, r_0/2)$ and all $Y \in \Omega^m$ such that
$\| X-Y\|\leq r_0/2$, we denote by
\[
\tilde{ \mathcal E}_{\sigma,Y} : {\mathcal C}^{0, \alpha} (\bar
\Omega_\sigma
 \, (Y)) \longrightarrow {\mathcal C}^{0, \alpha}_\nu (\bar \Omega^* \, (Y)),
\]
the extension operator defined by $\tilde {\mathcal E}_{\sigma,Y}
(f) = f$ in $\bar \Omega_\sigma \, (Y)$
\[
\tilde {\mathcal E}_{\sigma,Y} (f) \, (y^i+x) = \tilde \chi \left(
\frac{|x|}{\sigma} \right) \, f \left(y^i + \sigma \, \frac{x}{|x|}
\right)
\]
for each $j=1, \ldots, m$ and $\tilde {\mathcal E}_{\sigma,Y} (f) =
0$ in each $B_{\sigma/2} (y^j)$, where $t \longmapsto \tilde \chi
(t)$ is a cutoff function identically equal to $1$ for $t \geq 1$
and identically equal to $0$ for $t \leq 1/2$. It is easy to check
that there exists a constant $c =  c(\nu) >0$ only depending on
$\nu$ such that
\begin{equation}
\| \tilde {\mathcal E}_{\sigma,Y} \, (w) \|_{{\mathcal C}^{0,
\alpha}_{\nu} (\bar \Omega^* \, (X))} \leq c \, \| w \|_{{\mathcal
C}^{0, \alpha}_{\nu} (\bar \Omega_{\sigma} \, (X))}. \label{le:7.2}
\end{equation}
\medskip

We fix
\[
\nu \in (-1,0),
\]
and denote by ${\mathcal  G}_{\nu, Y}$ the right inverse provided by
Proposition~\ref{pr:4.1}. Clearly, it is enough to find $\tilde v
\in {\mathcal C}^{4, \alpha}_\nu (\Omega^* \, (Y))$ solution of
\begin{equation}
\tilde v =  \tilde N (\e, \Lambda , Y, \Phi, \Psi \, ; \, \tilde v )
\label{eq:7.30}
\end{equation} where
we have defined
\[
\tilde N (\e, \Lambda , Y, \Phi, \Psi \, ; \, \tilde v ) : =
{\mathcal G}_{\nu, Y} \circ {\tilde {\mathcal E}}_{r_\e , Y} \,
\left( \rho^{4} \, e^{{\bf \tilde u} + \tilde{v}} - \Delta^{2} \,
{\bf \tilde u} \right).
\]

Given $\kappa >0$ (whose value will be fixed later on), we further
assume that $\Phi$ and $\Psi$ satisfy
\begin{equation}
\| \Phi \|_{({\cal C}^{4, \alpha} (S^3))^m}\leq \kappa \, \e , \quad
\mbox{and} \quad \| \Psi \|_{({\cal C}^{2, \alpha} (S^3))^m }\leq
\kappa \, \e. \label{eq:7.4}
\end{equation}
Moreover, we assume that the parameters $\Lambda$ and the points $Y$
are chosen to satisfy
\begin{equation}
| \Lambda |\leq \kappa \, \varepsilon , \qquad \mbox{and} \qquad
\|Y-X \| \leq \kappa \,\sqrt \e . \label{eq:7.5}
\end{equation}

Then, the following result holds~:
\begin{lemma}
Given $\kappa >1$. There exist $\e_\kappa >0$, $c_\kappa
>0$ and $\bar c_\kappa >0$ such that, for all $\e \in (0,
\e_\kappa)$
\begin{equation}
\| \tilde N (\e, \Lambda , Y, \Phi, \Psi \, ; \, 0 ) \|_{{\mathcal
C}^{4, \alpha}_\nu (\bar \Omega^* \, (Y))}\leq c_\kappa \,
\varepsilon^{3/2}. \label{eq:7.6}
\end{equation}
Moreover,
\begin{equation}
\| \tilde N ( \e, \Lambda , Y, \Phi, \Psi \, ; \, \tilde v_2 ) -
\tilde N (\e, \Lambda , Y, \Phi, \Psi \, ; \, \tilde v_1 )
\|_{{\mathcal C}^{4, \alpha}_\nu (\bar \Omega^* \, (Y))}\leq \bar
c_\kappa \, \varepsilon^{2} \, \| \tilde v_2 - \tilde v_1
\|_{{\mathcal C}^{4, \alpha}_\nu (\bar \Omega^* \, (Y))}
\label{eq:7.7}
\end{equation}
and
\begin{equation}
\begin{array}{rllllll}
\| \tilde N (\e, \Lambda , Y, \Phi_2 , \Psi_2 \, ; \, \tilde v ) -
\tilde N (\e, \Lambda , Y, \Phi_1 , \Psi_1 \, ; \, \tilde v )
\|_{{\mathcal C}^{4, \alpha}_\nu (\bar \Omega^* \, (Y))} \qquad \qquad \qquad \qquad \\[3mm]
\leq \bar c_\kappa \, \varepsilon^{1/2} \, \left( \| \Phi_2 -
\Phi_1 \|_{({\mathcal C}^{4, \alpha}(S^3))^m} + \| \Psi_2 - \Psi_1
\|_{({\mathcal C}^{2, \alpha}(S^3))^m} \right)
\end{array}
\label{eq:7.8}
\end{equation}
provided $\tilde v= v, v_1, v_2 \in {\mathcal C}^{4, \alpha}_\nu
(\bar \Omega^* \, (Y))$, $\tilde \Phi = \Phi, \Phi_1, \Phi_2 \in
({\mathcal C}^{4, \alpha} (S^3))^m$, $\tilde \Psi = \Psi, \Psi_1,
\Psi_2 \in ({\mathcal C}^{2, \alpha} (S^3))^m$ satisfy
\[
\| \tilde v \|_{{\mathcal C}^{4, \alpha}_\nu ({\bar \Omega^*} \,
(Y))} \leq 2 \, c_\kappa \, \e^{3/2} , \qquad \|\tilde \Phi
\|_{({\cal C}^{4, \alpha} (S^3))^m}\leq \kappa \, \varepsilon ,
\qquad \| \tilde \Psi \|_{({\cal C}^{2, \alpha} (S^3))^m}\leq \kappa
\, \varepsilon ,
\]
and $| \Lambda | \leq \kappa \, \e $,  $\|Y -X\| \leq \kappa \,
\sqrt \e $. \label{le:7.1}
\end{lemma}
{\bf Proof~:} The proof of the first estimate follows from the
result of Lemma~\ref{le:5.2} together with (\ref{eq:7.4}). More
precisely, we have
\[ \|  \rho^{4} \, e^{{\bf \tilde u} } \|_{{\mathcal C}^{0, \alpha}_{\nu -4} (\bar
\Omega_{r_\e} \, (Y) ) } \leq c_\kappa \, \e^{(4- \nu)/2} \qquad
\mbox{and} \qquad \| \Delta^{2} \, {\bf \tilde u} \|_{{\mathcal
C}^{0, \alpha}_{\nu -4} (\bar \Omega_{r_\e} \, (Y) ) } \leq c_\kappa
\, \e^{3/2} .
\]
The proof of the first estimate follows from (\ref{eq:7.2}) and
Proposition~\ref{pr:4.1}.

\medskip

The proof of the second estimate follows from
\[ \|  \rho^{4} \,  ( e^{{\bf \tilde u} + v_2} -  e^{{\bf \tilde u} + v_1} )
\|_{{\mathcal C}^{0, \alpha}_{\nu -4} (\bar
\Omega_{r_\e} \, (Y) ) } \leq c_\kappa \, \e^{2} \, \| v_2-v_1
\|_{{\mathcal C}^{4, \alpha}_\nu (\bar \Omega^* \, (Y))}
\]
and the third estimate follows from
\[ \|  \rho^{4} \, ( e^{{\bf \tilde u}_2 + v} -  e^{{\bf \tilde u}_1 + v})
\|_{{\mathcal C}^{0, \alpha}_{\nu -4} (\bar \Omega_{r_\e} \, (Y) )
} \leq c_\kappa \, \e^{(4- \nu)/2} \, \left( \| \Phi_2 - \Phi_1
\|_{({\mathcal C}^{4, \alpha}(S^3))^m} + \| \Psi_2 - \Psi_1
\|_{({\mathcal C}^{2, \alpha}(S^3))^m} \right)
\]
and
\[
\| \Delta^{2} \, ({\bf \tilde u}_2 - {\bf \tilde u}_1)
\|_{{\mathcal C}^{0, \alpha}_{\nu -4} (\bar \Omega_{r_\e} \, (Y) )
} \leq c_\kappa \, \e^{1/2} \, \left( \| \Phi_2 - \Phi_1
\|_{({\mathcal C}^{4, \alpha}(S^3))^m} + \| \Psi_2 - \Psi_1
\|_{({\mathcal C}^{2, \alpha}(S^3))^m} \right)
\]
(where ${\bf \tilde u}_j$ corresponds to ${\bf \tilde u}$ when $\Phi
=\Phi_j$ and $\Psi = \Psi_j$) together with (\ref{eq:7.2}) and
Proposition~\ref{pr:4.1}. \hfill $\Box$

\medskip

Reducing $\e_\kappa$ is necessary, we can assume that \[ \bar
c_\kappa \, \e^2 \leq \frac{1}{2} \] for all $\e \in (0,
\e_\kappa)$. Then, (\ref{eq:7.6}) and (\ref{eq:7.7}) are enough to
show that \[ \tilde v \, \longmapsto \tilde N (\e, \Lambda , Y,
\Phi, \Psi \,  \, \tilde v)
\]
is a contraction from
\[
\{ v \in {\mathcal C}^{4, \alpha}_\nu (\bar \Omega^* \, (Y))  \qquad
: \qquad \|v \|_{{\mathcal C}^{4, \alpha}_\nu (\bar \Omega^* \,
(Y))} \leq 2 \, c_\kappa \, \e^{3/2} \}
\]
into itself and hence has a unique fixed point $\tilde v (\e,
\Lambda , Y, \Phi, \Psi \, ; \, \cdot )$ in this set. This fixed
point is {\underline a} solution of (\ref{eq:7.3}). Reducing
$\e_\kappa$ if this is necessary, it follows from (\ref{eq:7.7}) and
(\ref{eq:7.8}) in Lemma~\ref{le:7.1} that
\begin{equation}
\begin{array}{rlllll}
\| \tilde v (\e, \Lambda , Y, \Phi_2 , \Psi_2 \, ; \, \cdot ) -
\tilde v (\e, \Lambda , Y, \Phi_1, \Psi_1 \, ; \, \cdot ) \| \qquad \qquad \qquad \\[3mm]
\leq 2 \, \bar c_\kappa \, \varepsilon^{1/2} \, \left( \| \Phi_2 -
\Phi_1 \|_{({\mathcal C}^{4, \alpha}(S^3))^m} + \| \Psi_2 - \Psi_1
\|_{({\mathcal C}^{2, \alpha}(S^3))^m} \right). \end{array}
\label{eq:7.9}
\end{equation}

We summarize this in the~:
\begin{proposition}
Given $\kappa > 0$, there exists $\e_\kappa >0$ and $c_\kappa >0$
(only depending on $\kappa$) such that for all $\e \in (0,
\e_\kappa)$, for all set of parameters $\Lambda$, points $Y$
satisfying
\[
| \Lambda |\leq \kappa \, \varepsilon , \qquad \mbox{and} \qquad
\|Y-X \| \leq \kappa \,\sqrt \e
\]
and boundary functions $\Phi$ and $\Psi$ satisfying (\ref{eq:5.6})
and
\[
\| \Phi \|_{({\cal C}^{4, \alpha} (S^3))^m}\leq \kappa \, \e , \quad
\mbox{and} \quad \| \Psi \|_{({\cal C}^{2, \alpha} (S^3))^m }\leq
\kappa \, \e.
\]
the function
\[
\begin{array}{rllllll}
\tilde u(\e, \Lambda , Y, \Phi, \Psi \, ; \, \cdot )  & : = &
\displaystyle \sum_{j=1}^m (1 + \lambda^j ) \, G_{y^j} +
\sum_{j=1}^m \chi_{r_0} (\cdot - y^j) \, H^{e} (\varphi^j, \psi^j \,
; \, (\cdot - y^j)/
r_\varepsilon ) \\[3mm]
 &  + & \tilde v(\e, \Lambda , Y, \Phi, \Psi \, ; \,
\cdot ),
\end{array}
\]
solves (\ref{eq:7.2}) in $\bar \Omega_{r_\e} \, (Y)$. In addition
\begin{equation}
\| \tilde v (\e, \Lambda , Y, \Phi, \Psi \, ; \, \cdot ) \|_{{\cal
C}^{4,\alpha}_{\nu}(\bar \Omega^*)} \leq 2 \, c_\kappa \, \e^{3/2}.
\label{eq:7.10} \end{equation}
\label{pr:7.1}
\end{proposition}

Observe that the function $ \tilde v_{\e, \Lambda , Y , \Phi, \Psi}$
being obtained as a fixed point for contraction mapping, it depends
continuously on the parameters $\Lambda $ and the points $Y$.

\section{The nonlinear Cauchy-data matching}

Keeping the notations of the previous sections, we gather the
results of the Proposition~\ref{pr:6.1} and
Proposition~\ref{pr:7.1}. From now on $\kappa >1$ is fixed large
enough (we will shortly see how) and $\e \in (0, \e_\kappa)$.

\medskip

Assume that $X =(x^1, \ldots, x^m) \in \Omega^m$  is a nondegenerate
critical point of the function $W$ defined in the introduction. For
all $j=1, \ldots, m$, we define $\tau^j_*
>0$ by
\begin{equation}
- 4 \, \log \tau^j_* =  R (x^j, x^j) + \sum_{\ell \neq j} G(x^\ell,
x^j). \label{eq:8.0}
\end{equation}

We assume that we are given~:
\begin{itemize}

\item[(i)] points $Y : = (y^1, \ldots, y^m) \in \Omega^m$ close to $X : = (x^1,
\ldots, x^m)$ satisfying (\ref{eq:7.5}).

\item[(ii)] parameters $\Lambda : = (\lambda^1, \ldots, \lambda^m) \in {\mathbb R}^m
$ satisfying (\ref{eq:7.5}).

\item[(iii)] parameters $T : = (\tau^1, \ldots, \tau^m) \in (0,
\infty)^m$ satisfying (\ref{eq:6.5}) (where, for each $j=1, \ldots,
m$, $\tau_*$ is replaced by $\tau_*^j$).

\end{itemize}

We set
\[
R_\e^j: = \tau^j/ \sqrt \e
\]
First, we consider some set of boundary data
\[
\Phi : = (\varphi^1, \ldots, \varphi^m) \in ({\cal
C}^{4,\alpha}(S^3))^m \qquad \mbox{and} \qquad \Psi : = (\psi^1,
\ldots, \psi^m) \in ({\cal C}^{2,\alpha}(S^3))^m
\]
satisfying (\ref{eq:5.2}) and (\ref{eq:6.5}).

\medskip

Thanks to the result of Proposition~\ref{pr:6.1}, we can find
$u_{int}$ a solution of
\[
\Delta^2 u - \rho^4 \, e^u =0
\]
in each $B_{r_\e} (y^j)$, which can be decomposed as
\[
\begin{array}{rlllll}
u_{int} (\e, T,Y, \Phi, \Psi \, ; \, x) & : = &
u_{\varepsilon,\tau^j}(x- y^j) + H^{i} ({ \varphi^j}, { \psi^j} \, ;
\, (x-y^j)/r_\varepsilon ) \\[3mm]
& + & v (\e, \tau^j, \varphi^j, \psi^j \, ; \, R_\e^j (x - y^j)/
r_\e )
\end{array}
\]
in $B_{r_\e} (y^j)$.

\medskip

Similarly, given some boundary data
\[
\tilde \Phi  : = (\tilde \varphi^1, \ldots , \tilde \varphi^m ) \in
({\cal C}^{4,\alpha}(S^3))^m \qquad \mbox{and} \qquad \tilde \Psi :=
(\tilde \psi^1, \ldots, \tilde \psi^m ) \in ({\cal
C}^{2,\alpha}(S^3))^m
\]
satisfying (\ref{eq:5.6}) and (\ref{eq:7.4}), we use the result of
Proposition~\ref{pr:7.1}, to find $u_{ext}$ a solution of
\[
\Delta^2 u - \rho^4 \, e^u =0
\]
in $\bar \Omega_{r_\e} \, (Y)$, which can be decomposed as
\[
\begin{array}{rllllll}
u_{ext}(\e, \Lambda, Y, \tilde \Phi, \tilde \Psi \, ;  \, x) & =&
\displaystyle \sum_{j=1}^m (1+ \lambda^j) \, G (y^j , x) +
\sum_{j=1}^m \chi_{r_0} (x - y^j) \,  H^{e} (\tilde{{\varphi}}^j, {
\tilde{\psi}}^j \, ; \, (x - y^j) / r_\varepsilon ) \\[3mm]
 & + & {\tilde v} (\e, \Lambda, Y, \tilde \Phi,
\tilde \Psi \, ; \, x). \end{array}
\]

It remains to determine the parameters and the boundary functions in
such a way that the function which is equal to $u_{int }$ in $\cup_j
\, B_{r_\e} (y^j)$ and which is equal to $u_{ext}$ in $\bar
\Omega_{r_\e} \, (Y)$ is a smooth function. This amounts to find the
boundary data and the parameters so that, for each $j=1, \ldots, m$
\begin{equation}
u_{int} = u_{ext} , \qquad \partial_r u_{int} =
\partial_r u_{ext}, \qquad \Delta u_{int} = \Delta u_{ext} , \qquad \partial_r \Delta u_{int} =
\partial_r \Delta u_{ext},
\label{eq:8.1}
\end{equation}
on $\del B_{r_\e} (y^j)$. Assuming we have already done so, this
provides for each $\e$ small enough a function $u_\e \in {\mathcal
C}^{4, \alpha}(\bar \Omega)$ (which is obtained by patching together
the function $u_{int}$ and the function $u_{ext}$) solution of
$\Delta^2 \, u  - \rho^4 \, e^u = 0$ and elliptic regularity theory
implies that this solution is in fact smooth. This will complete the
proof of our result since, as $\e$ tends to $0$, the sequence of
solutions we have obtained satisfies the required properties,
namely, away from the points $x^j$ the sequence $u_\e$ converges to
$\sum_j G (x^j , \, \cdot \, )$.

\medskip

Before, we proceed, some remarks are due. First it will be
convenient to observe that the functions $u_{\e, \tau^j}$ can be
expanded as
\begin{equation}
u_{\e, \tau^j} (x) = - 8 \, \log |x| - 4 \log \tau^j + {\mathcal
O}(\e) \label{eq:8.2}
\end{equation}
 near $\del B_{r_\e}$.
Also, the function
\[
\sum_{j=1}^m (1+ \lambda^j) \, G (y^j  ,x)
\]
which appears in the expression of $u_{ext}$ can be expanded as
\begin{equation} \sum_{\ell=1}^m (1+ \lambda^\ell) \,
G (y^\ell, y^j+x ) = - 8 \, (1+ \lambda^j) \, \log |x| + E_j (Y ;
y^j) + \nabla E_j (Y ; y^j) \cdot x + {\mathcal O} (\e)
\label{eq:8.3}
\end{equation} near $\del B_{r_\e} (y^j)$. Here, we have defined
\[
E_j(Y ; \cdot ) : = R(y^j,  \, \cdot  )+ \sum_{\ell \neq j} G
(y^\ell , \, \cdot).
\]

In (\ref{eq:8.1}), all functions are defined on $\del B_{r_\e}
(y^j)$, nevertheless, it will be convenient to solve, instead of
(\ref{eq:8.1}) the following set of equations
\begin{equation}
\begin{array}{rrrrrrrllllllll}
(u_{int} -u_{ext} )(y^j + r_\e \, \cdot )  & = & 0, & \qquad &
(\partial_r u_{int } - \partial_r u_{ext}) (y^j + r_\e \, \cdot ) & = & 0, \\[3mm]
(\Delta  u_{int } - \Delta u_{ext}) (y^j + r_\e \, \cdot ) & = & 0 ,
& \qquad &  ( \partial_r \, \Delta u_{int } - \partial_r \, \Delta
u_{ext} ) (y^j + r_\e \, \cdot ) & = & 0,
\end{array}
\label{eq:8.4}
\end{equation}
on $S^3$. Here all functions are considered as functions of $z \in
S^3$ and we have simply used the change of variables $x = y^j + r_\e
\, z$ to parameterize $\del B_{r_\e} (y^j)$.

\medskip

Since the boundary data satisfy (\ref{eq:5.2}) and (\ref{eq:5.6}),
we decompose
\[
\Phi  = \Phi_0 + \Phi_1 + \Phi^{\perp} \qquad \qquad \Psi  = 8 \,
\Phi_0 + 12 \, \Phi_1 + \Psi^{\perp}
\]
and
\[
\tilde \Phi = \tilde \Phi_0 + \tilde \Phi_1 + \tilde \Phi^{\perp}
\qquad \qquad \tilde \Psi = \tilde \Psi_1 + \tilde \Psi^{\perp}
\]
where the components of $\Phi_0, \tilde \Phi_0$ are constant
functions on $S^3$, the components of $\Phi_1 , \tilde \Phi_1,
\tilde \Psi_1$ belong to $\mbox{Ker} (\Delta_{S^3} +3) = \mbox{Span}
\{e_1, \ldots, e_4\}$ and where the components of $\Phi^{\perp} ,
\Psi^{\perp}, \tilde \Phi^{\perp} , \tilde \Psi^{\perp}$ are $L^2
(S^3)$ orthogonal to the constant function and the functions $e_1,
\ldots, e_4$. Observe that the components of $\Psi$ over the
constant functions or functions in $\mbox{Ker} (\Delta_{S^3} +3)$
are determined by the corresponding components of $\Phi$. Moreover,
$\tilde \Psi$ has no component over constant functions.

\medskip

We first consider the $L^2(S^3)$-orthogonal projection of
(\ref{eq:8.4}) onto the space of functions which are orthogonal to
the constant function and the functions $e_1, \ldots, e_4$. This
yields the system
\begin{equation}
\left\{ \begin{array}{rllllll}
\varphi^{j,\perp} -\tilde \varphi^{j,\perp} & = &  M_0^{(j)} (\e, \Lambda, T, Y, \Phi, \tilde \Phi, \Psi, \tilde \Psi) \\[3mm]
\del_r \, H^i (\varphi^{j,\perp} , \psi^{j,\perp} \, ; \cdot) -
\del_r \,  H^e (\tilde \varphi^{j,\perp} , \tilde \psi^{j,\perp} \, ; \cdot) & = & M_1^{(j)} (\e, \Lambda, T, Y, \Phi, \tilde \Phi, \Psi, \tilde \Psi )   \\[3mm]
\psi^{j,\perp} - \tilde \psi^{j,\perp} &= & M_2^{(j)} (\e, \Lambda, T, Y, \Phi, \tilde \Phi, \Psi, \tilde \Psi ) \\[3mm]
\del_r \, \Delta \, H^i (\varphi^{j,\perp}, \psi^{j,\perp} \, ;
\cdot) - \del_r \, \Delta \, H^e ( \tilde \varphi^{j,\perp} , \tilde \psi^{j,\perp} \, ; \cdot) ) & = & M_3^{(j)} (\e, \Lambda, T, Y, \Phi, \tilde \Phi, \Psi, \tilde \Psi ) \\[3mm]
\end{array}
\right. \label{eq:8.5}
\end{equation}
where the functions $M_k^{(j)}$ are nonlinear functions of the
parameters $\e$, $\Lambda$, $Y$, $T$ and the boundary data $\Phi$,
$\tilde \Phi$, $\Psi$ and $\tilde \Psi$. Moreover, using
(\ref{eq:8.2}) and (\ref{eq:8.3}) and also (\ref{eq:6.11}) (keeping
in mind that $\mu \in (1,2)$) and (\ref{eq:7.10}) (keeping in mind
that $\nu \in (-1,0)$), we conclude that, for each $j=1, \ldots, m$
and $k = 0, 1, 2, 3$
\begin{equation}
\|M_k^{(j)}\|_{{\mathcal C}^{4-k, \alpha} (S^3)} \leq c\, \e
\label{eq:8.55}
\end{equation}
for some constant $c >0$ independent of $\kappa$ (provided $\e \in
(0, \e_\kappa)$).

\medskip

Thanks to the result of Lemma~\ref{le:5.3} and (\ref{eq:8.55}), this
last system can be re-written as
\[
(\Phi^{\perp} , \tilde \Phi^{\perp} , \Psi^{\perp}, \tilde
\Psi^{\perp} ) = M(\e, \Lambda, T, Y, \Phi, \tilde \Phi, \Psi,
\tilde \Psi )
\]
where
\[
\|M \|_{({\mathcal C}^{4, \alpha} (S^3) )^{2m} \times ({\mathcal
C}^{2, \alpha} (S^3) )^{2m}} \leq c\, \e
\]for some
constant $c >0$ independent of $\kappa$ (provided $\e \in (0,
\e_\kappa)$). Moreover, (\ref{eq:6.10}) and (\ref{eq:7.9}) imply
(reducing $\e_\kappa$ if necessary) that, the mapping $M$ is a
contraction from the ball of radius $\kappa \, \e$ in $({\mathcal
C}^{4, \alpha} (S^3) )^{2m} \times ({\mathcal C}^{2, \alpha} (S^3)
)^{2m}$ into itself and as such has a unique fixed point in this
set. Observe that this fixed point depends continuously on $\e$,
$\Lambda$, $T$, $Y$ and also on $\Phi_0$, $\tilde \Phi_0$, $\Phi_1$,
$\tilde \Phi_1$ and $\tilde \Psi_1$.

\medskip

We insert this  fixed point in (\ref{eq:8.4}) and now project the
corresponding system over the set of functions spanned by $e_1,
\ldots, e_4$ and finally over the set of constant functions.

\medskip

The first projection yields the system of equations
\begin{equation}
\left\{ \begin{array}{rlllll}
\Phi_1 & = & \bar M_{1}(\e, \Lambda, T, Y,\Phi_0, \tilde \Phi_0,  \Phi_1, \tilde \Phi_1, \tilde \Psi_1 )  \\[3mm]
\tilde \Phi_1 & = & \bar M_2 (\e, \Lambda, T, Y , \Phi_0, \tilde \Phi_0,  \Phi_1, \tilde \Phi_1, \tilde \Psi_1)  \\[3mm]
\Psi_1  & = & \bar M_3 (\e, \Lambda, T, Y,\Phi_0, \tilde \Phi_0,  \Phi_1, \tilde \Phi_1, \tilde \Psi_1 )  \\[3mm]
\sqrt\e \, \nabla E_j (Y ; y^j)  & = & \bar M_4^{(j)} (\e, \Lambda,
T, Y ,\Phi_0, \tilde \Phi_0,  \Phi_1, \tilde \Phi_1, \tilde \Psi_1)
\end{array}
\right. \label{eq:8.6}
\end{equation}
where the functions $\bar M_k$ (and also $\bar M_4^{(j)}$) are
nonlinear functions depending continuously on the parameters $\e$,
$\Lambda$, $T$, $Y$ and the components of the boundary data
$\Phi_0$, $\tilde \Phi_0$, $\Phi_1$, $\tilde \Phi_1$ and $\tilde
\Psi_1$. Moreover,
\[
|\bar M_k | \leq c\, \e
\]
for some constant $c >0$ independent of $\kappa$ (provided $\e \in
(0, \e_\kappa)$).

\medskip

Let us comment briefly on how these equations are obtained. These
equations simply come from (\ref{eq:8.1}) when expansions
(\ref{eq:8.2}) and (\ref{eq:8.3}) are taken into account, together
with the expression of $H^i (\varphi^j, \psi^j  \, ; \cdot)$ and
$H^{e} (\tilde \varphi^j , \tilde \psi^j \, ; \cdot)$ given in
Lemma~\ref{le:5.1} and Lemma~\ref{le:5.2}, and also the estimates
(\ref{eq:6.11}) and (\ref{eq:7.10}). Observe that the projection of
the term $x \longrightarrow \nabla E_j (Y ; y^j) \cdot x$ which
arises in (\ref{eq:8.3}), as well as the projection of its partial
derivative with respect to $r$, over the set of constant function is
equal to $0$. Moreover, this term projects identically over the set
of functions spanned by $e_1, \ldots, e_4$ as well as its derivative
with respect to $r$. Finally, its Laplacian vanishes identically.

\medskip

Recall that we have define in the introduction the function
\[
W (Y) : =  \sum_{j=1}^m R \, (y^j, y^j )+ \sum_{j_1\neq j_2} G
(y^{j_1}, y^{j_2})
\]
Using the symmetries of the functions $G$ and $R$, namely the fact
that
\[
G(x,y) = G(y,x) \qquad \mbox{and} \qquad R(x,y) = R(y,x)
\]
we get
\[
\nabla W|_Y = 2 \, (\nabla E_1 (Y ; y^1), \ldots, \nabla E_m (Y,
y^m)).
\]
Now, we have assumed that the point  $X= (x^1, \ldots, x^m)$ is a
nondegenerate critical point of the functional $W$ and hence
\[
\nabla W |_X =0,
\]
and
\[
({\mathbb R}^4)^m \ni Z \longmapsto D (\nabla W)|_{X}  (Z) \in (
{\mathbb R}^4)^m
\]
is invertible. Therefore, the last equation can be rewritten as
\[
\sqrt\e \, (Y-X) = \bar M_5 \, (\e, \Lambda, T, Y ,\Phi_0, \tilde
\Phi_0, \Phi_1, \tilde \Phi_1, \tilde \Psi_1)
\]

The projection of (\ref{eq:8.4}) over the constant function, leads
to the system
\begin{equation}
\left\{
\begin{array}{rlllll} (\log 1/ \e)^{-1} \, \log (\tau^j/ \tau^j_*) & = & \bar M_6 (\e, \Lambda, T, Y ,\Phi_0, \tilde \Phi_0,  \Phi_1, \tilde \Phi_1, \Psi_1, \tilde \Psi_1) \\[3mm]
\tilde \Phi_0 & = & \bar M_7 (\e, \Lambda, T, Y ,\Phi_0, \tilde \Phi_0,  \Phi_1, \tilde \Phi_1, \Psi_1, \tilde \Psi_1)  \\[3mm]
\Phi_0 & = & \bar M_8 (\e, \Lambda, T, Y ,\Phi_0, \tilde \Phi_0,  \Phi_1, \tilde \Phi_1, \Psi_1, \tilde \Psi_1) \\[3mm]
\Lambda & = & \bar M_9 (\e, \Lambda, T, Y ,\Phi_0, \tilde \Phi_0,
\Phi_1, \tilde \Phi_1, \Psi_1, \tilde \Psi_1)
\end{array}
\right. \label{eq:8.7}
\end{equation}
where the function $\bar M_k$ satisfy the usual properties. If we
define the parameters $U : = (u^1, \ldots, u^m)$ where
\[
u^j = \frac{1}{\log 1/\e} \, \log (\tau^j/\tau^j_*)
\]
and
\[
Z = \sqrt\e \, (Y-X)
\]
so that the system we have to solve
reads
\begin{equation}
( U , \Lambda , Z , \Phi_0 ,\tilde \Phi_0 , \Phi_1 , \tilde \Phi_1
, \tilde\Psi_1 ) = \bar M (\e, U , \Lambda , Z , \Phi_0 ,\tilde
\Phi_0 , \Phi_1 , \tilde \Phi_1 , \tilde \Psi_1 ). \label{eq:8.8}
\end{equation}
where as usual, the nonlinear function $\bar M$ depends continuously
on the parameters $T  , \Lambda , Z$ and the functions $\Phi_0 ,
\tilde \Phi_0 , \Phi_1 , \tilde \Psi_1 $ and is bounded (in the
appropriate norm) by a constant (independent of $\e$ and $\kappa$)
time $\e$, provided $\e \in (0, \e_\kappa)$. Observe that
\[
\begin{array}{crlllll}
U , \Lambda \in  {\mathbb R}^m , \qquad  Z \in ({\mathbb R}^4)^m ,
\qquad \Phi_0 , \tilde \Phi_0 \in {\mathbb R}^m
\\[3mm]
\Phi_1 ,  \tilde \Phi_1 , \tilde \Psi_1  \in (\mbox{Ker} \,
(\Delta_{S^3}+3))^m .
\end{array}
\]
In addition, reducing $\e_\kappa $ if necessary, this nonlinear
mapping sends the ball of radius $\kappa \, \e$ (for the natural
product norm) into itself, provided $\kappa$ is fixed large enough
and $\e \in (0, \e_\kappa)$. Applying Schauder's fixed point Theorem
in the ball of radius $\kappa \, \e$ in the product space where the
entries live yields the existence of a solution of (\ref{eq:8.8})
and this completes the proof of Theorem~\ref{th:1.1}.

\section{Comments}

Let us comment on how the condition "$(x^1,\ldots ,x^m)$ is a {\em
nondegenerate} critical point of $W$" enters in our analysis since,
we confess, that it is somehow very well hidden.

\medskip

The condition "$(x^1,\ldots ,x^m)$ is a critical point of $W$"
enters in the estimate (\ref{eq:8.3}) when $Y=X$ and $\Lambda =0$,
since, in this case we have
\[
 \sum_{\ell=1}^m G (x^\ell, x^j+x ) =
- 8 \, \log |x| + E_j (X ; x^j) + {\mathcal O} (\e)
\]
while, if $(x^1,\ldots ,x^m)$ were not a critical point of $W$, then
$\nabla E_j (X; x^j) \neq 0$ and we would only have
\[
\sum_{\ell=1}^m G (x^\ell, x^j+x ) = - 8 \, \log |x| + E_j (X ; x^j)
+ {\mathcal O} (\e^{1/2})
\]
which would not be enough : roughly speaking this says that the
approximate solution we have constructed is not close to any
solution of the problem. Given the result of Lin and Wei
\cite{Lin-Wei}, the condition on "$(x^1,\ldots ,x^m)$ being a
critical point of $W$" is a natural one.

\medskip

The origin of the "nondegeneracy" assumption takes its roots in the
result of Lemma~\ref{le:3.1} which classifies all the solutions of
the linearized equation about the rotationally symmetric solution.
The existence of elements $\phi_i$, for $i=1, \ldots, 4$ in the
kernel of ${\mathbb L}$ has forced us in proposition~\ref{pr:3.1} to
work with weights $\mu >1$ to obtain the surjectivity of the
operator $L_\mu$. This choice has one importance consequence : In
Lemma~\ref{le:5.1}, we had to restrict our attention to boundary
data which satisfy the constraints (\ref{eq:5.2}) and (\ref{eq:5.6})
(even though only the second constraint in (\ref{eq:5.2}) is
important to understand where the nondegeneracy condition comes
from) to obtain bi-harmonic extensions in the unit ball which vanish
at the origin at least quadratically. A second reading will convince
the reader that this property was crucial in the estimate of
Lemma~\ref{le:6.1}. Indeed, the main estimate in this Lemma arises
from the fact that
\[|H^i(\varphi, \psi ; \cdot /R_\e) | \leq c_\kappa^{(1)} \, \e^2 \,
|x|^2. \]
 Without the second hypothesis in
(\ref{eq:5.2}) we would only have \[ |H^i(\varphi, \psi ; \cdot
/R_\e) | \leq c_\kappa^{(1)} \, \e^{3/2} \, |x|\] which would have
led  in Lemma~\ref{le:6.1} to the estimate \[\| N (\e, \tau,
\varphi, \psi ; 0)\|_{{\mathcal C}^{4, \alpha}_\mu ({\mathbb R})^4}
\leq c_\kappa \, \e^{3/2}\] But since $\mu \in (1,2)$ this implies
that, on the boundary $\partial B_{R_\e}$ the function $v(\e, \tau,
\varphi, \psi ; \, \cdot \, )$ is bounded by a constant times
$\e^{(3-\mu)/2}$ and since
\[
\e^{(3-\mu)/2} >> \e
\]
the function $v$ would be much larger than the functions
$H^i(\varphi, \psi ; \, \cdot /R_\e)$ on this boundary and hence
could not be considered as a small perturbation anymore. Given the
fact that, in the construction of $H^i$ and $H^e$ we could not
prescribe any function, we had to "find" new degrees of freedom to
compensate the constraints imposed by (\ref{eq:5.2}) and
(\ref{eq:5.6}). The introduction of the parameters $\tau^j $ and
$\lambda^j$ enter at this point to overcome the fist condition
imposed by (\ref{eq:5.2}) and also the condition imposed by
(\ref{eq:5.6}). The points $y^j$ close to $x^j$ are introduced to
compensate the second condition imposed by (\ref{eq:5.2}) and this
is precisely were the nondegeneracy of the critical points of $W$
comes into play.

\medskip

Let us point out that the {\em nondegeneracy} condition strictly
speaking can be weakened as this has been done for example in
\cite{del-Mus} and \cite{Esp-Gro-Pis} in the case of equation
(\ref{eq:1.8}). The idea being that the nondegeneracy is essentially
used to solve the last equation in (\ref{eq:8.6}) by some disguised
version of the Implicit Function Theorem. But, remembering that the
problem we want to solve is a variational problem, this last
equation can be rephrased essentially as the gradient of a function
$W_\e$ which is defined on $\Omega^k$ and which converges (in a
sense to be made precise) to the function $W$ as $\e$ tends to $W$.
Some extra work is needed, but in any case, we could have used some
variational technics to find critical points of this functional.
Since nondegeneracy of critical points is a generic condition and in
order not to make the exposition of this "nonlinear domain
decomposition technic" as clear as possible, we have chosen not to
follow this route.

\section{Further results}

Modifying very little the previous analysis, it is possible to
extend the result of Theorem~\ref{th:1.1} to handle more general
equations. We will illustrate this on one example.

\medskip

As usual, let us assume that $\Omega \subset {\R^4}$ is a regular
bounded open subset and let us choose $z^1, \ldots , z^p \in \Omega$
and  $\alpha^1,\ldots, \alpha^p \in (0, +\infty)$. We would like to
extend the result of Theorem~\ref{th:1.1} to the equation
\begin{equation}
\left\{\begin{array}{rclll}\Delta^{2} u &=& \displaystyle \rho^{4}\,
e^{u} - 64 \, \pi^2 \, \sum_{i=1}^{p} \alpha^i \, \delta_{z^i}
&\mbox{in}&
\Omega\\[3mm]
u & = & \Delta u =0 & \mbox{on} & \partial\Omega.
\end{array}
\right. \label{eq:9.1}
\end{equation}
Namely, we are still looking for solutions which concentrate at some
points $x^1 , \ldots , x^m \in \Omega$, as the parameter
$\rho\longrightarrow 0$ and, in order to keep the amount of
technicalities as low as possible, we will assume that the set of
concentration points $x^j$ and the set of singularities $z^i$ are
disjoint. This problem is very much in the spirit of the work of
\cite{Esp-Gro-Pis} and \cite{del-Mus} even though we do no know any
applications in physics. On the other end solutions of this problem
might be of interest to understand constant $Q$-curvature metrics
with conical singularities.

\medskip

Setting
\[
v : = u +  \frac{1}{2}\sum_{i=1}^{p} \alpha^i \, G(z^i, \, \cdot)
\]
we can rephrase the equation satisfied by $u$ as an equation
satisfied by $v$, namely
\begin{equation}\left\{\begin{array}{rclll} \Delta^{2} v & = & \displaystyle \rho^{4}
\prod_{i=1}^{p}|x-z^i|^{4 \, \alpha^i}  \, e^{v} &\mbox{in}&
\Omega\\[3mm]
v& =&\Delta v=0 & \mbox{on}& \partial\Omega .
\end{array}
\right. \label{eq:9.2}
\end{equation}

This equation is a particular case of the more general problem
\begin{equation}\left\{\begin{array}{rclll}\Delta^{2} u &=&
\rho^{4} \, V \, e^{u} &\mbox{in}& \Omega \\[3mm]
u & = & \Delta u = 0  & \mbox{on}& \partial\Omega,
\end{array}
\right. \label{eq:9.3}
\end{equation}
where $V : \Omega \longrightarrow [0, +\infty)$ is a smooth
function. We are still looking for solutions of this last equation
which concentrate at some points $x^1, \ldots , x^m$, as the
parameter $\rho\longrightarrow 0$. In order to keep the
technicalities as low as possible, we will assume that the set of
concentration points $x^j$ and the set of zeros of $V$ are disjoint.

\medskip

As in the introduction, we introduce the functional
\begin{equation}
W (x^1,\ldots, x^m) : = \sum_{j=1}^{m} \, R(x^j, x^j) + \sum_{j \neq
\ell} G(x^j ,x^\ell) + 2 \, \sum_{j=1}^{m} \log V (x^j).
\label{eq:9.4}
\end{equation}
It is easy to check that the result of Theorem~\ref{th:1.1} holds
when (\ref{eq:1.1}) is replaced by (\ref{eq:9.3}) and (\ref{eq:9.4})
replaces (\ref{eq:1.4}). We briefly describe the main modifications
which are needed to prove this modified result.

\medskip

Only Sections 6,7 and 8 have to be slightly modified. In Section 6,
(\ref{eq:6.1}) has to be replaced by
\[
\Delta^2 \, u =  24 \, e^u  + \e^4 \, g
\]
where $g$ is a bounded function (in fact bounded in ${\mathcal
C}^{0, \alpha}(B_{R_\e})$ by some constant independent of $\e$). It
is easy to check that the analysis goes through. The presence of the
term $\e^4 \, g$ does not alter the estimates of Lemma~\ref{le:6.1}
and in fact, keeping the notations of introduced in the proof of
Lemma~\ref{le:6.1}, we have
\[
\| \e^4 \, g\|_{{\mathcal C}^{0, \alpha}_{\mu-4}(B_{R_\e})} \leq c
\, \e^{2 + \mu/2}.
\]
The result of Proposition~\ref{pr:6.1} remains unchanged. Section 7
applies {\it vertabim} and Proposition~\ref{pr:7.1} is unchanged.

\medskip

In Section 8, the main modification due is in the definition of
$u_{int}$. Indeed, for each $j=1, \ldots, m$ we apply the result of
the modified version of Section 6 with
\[
g = \frac{1}{\tau_j^4} \, (\Delta^2 \, \log V) (y^j + \e \, \cdot
/\tau_j).
\]
This induces in each $B_{r_\e} (y^j)$ a solution of
\[
\Delta^2 u = \rho^4 \,V \, e^u
\]
which can be decomposed as
\[
u_{int}(x) = u_{\varepsilon,\tau^j}(x-y^j) - \log V (x) + H^{i}(
\varphi^j, \psi^j\, ; \, (x -y^j)/r_\varepsilon ) + v (\e, \tau^j,
\varphi^j, \psi^j\, ; \, R_\e^j (x - y^j)/ r_\e) .
\]
The remaining of the analysis of Section 8 remains essentially
unchanged once the definition of $E_j$ is modified into
\[
E_j(Y ; \cdot ) : = R(y^j,  \, \cdot  )+ \sum_{\ell \neq j} G
(y^\ell , \, \cdot)  + \log V(y^j).
\]
We leave the details to the reader.

Sami Baraket. D\'epartement de Math\'ematiques, Facult\'e des
Sciences de Tunis, Tunisie. \\

Makkia Dammak. D\'epartement de Math\'ematiques, Facult\'e des
Sciences de Tunis, Tunisie.\\

Taieb Ouni. D\'epartement de Math\'ematiques, Facult\'e des Sciences
de Tunis, Tunisie.\\

Frank Pacard. Universit\'e Paris 12 et Institut Universitaire de
France.\\

\end{document}